%% file: paperANOR.tex
\RequirePackage{fix-cm}
\documentclass[smallextended]{svjour3}       
\smartqed  
\usepackage{graphicx}
\usepackage{subfigure}
\usepackage{ulem}
\usepackage{xcolor}
\usepackage{algorithm}
\usepackage{algpseudocode}
\usepackage{ulem}
\usepackage[urlcolor=black]{hyperref}

\newcommand{\Ii}{\mathbf{1}}
\begin{document}
\title{Inferences
for 
Random Graphs Evolved by Clustering Attachment
}


\author{Natalia Markovich \and Maksim Ryzhov \and  Marijus Vai{\v c}iulis
}


\institute{N. Markovich \at
              V.A.Trapeznikov Institute of Control Sciences,  Russian Academy of Sciences \\
              Tel.: +7(495)3348820\\
              \email{nat.markovich@gmail.com}           
\and
           M. Ryzhov \at
              V.A.Trapeznikov Institute of Control Sciences,  Russian Academy of Sciences \\
              Tel.: +7(495)3348820\\
              \email{ryzhov@phystech.edu}
              \and
           M. Vai{\v c}iulis \at
              Vilnius University, Institute of Data Science and Digital Technologies\\ Akademijos st. 4, LT-08663 Vilnius, Lithuania\\ \email{marijus.vaiciulis@mif.vu.lt}
}

\date{Received: date / Accepted: date}

\maketitle

\begin{abstract}
The evolution of random undirected graphs by the clustering attachment (CA) both  without node and edge deletion and with uniform node or edge deletion is investigated. Theoretical results  are obtained for the CA  without  node and edge deletion
when a newly appended node is connected to two existing nodes of the graph at each evolution step.
Theoretical results concern to 
(1)  the sequence of increments of the consecutive mean clustering coefficients tends to zero;
(2)  the sequences of node degrees and triangle counts of any fixed node which are proved to be submartingales. These results were obtained for any initial graph. The simulation study is provided  for the CA   with  uniform node or edge deletion and without any deletion. It is shown that
(1) the CA leads to light-tailed distributed node degrees and triangle counts;
(2) the average clustering coefficient  tends   to a constant
over time;
(3) the mean node degree and the mean triangle count increase over time with the rate depending on the parameters of the CA. The exposition is accompanied by a real data study.
\end{abstract}
\input{Sec1}
\input{Sec2ANOR}
\input{Sec3ANOR}
\input{Sec4ANOR}
\input{Sec5ANOR}
\input{Sec6}
\input{appendix.tex}

\input{reference.tex}
\begin{acknowledgements}
N. Markovich and M. Ryzhov were supported by the Russian Science Foundation
(grant \mbox{No.\,24-21-00183}). 
\end{acknowledgements}
\end{document}

%% file: Sec1.tex
\section{Introduction}\label{Sec1}
\par 
Network evolution attracts interest of researchers due to numerous applications (Avrachenkov \& Dreveton, 2022; Ghoshal et al., 2013; van der Hofstad, 2017; Norros \& Reittu, 2006; Wan et al., 2020).
The  popular mechanism to model  growing real-world networks and to explain a power-law distribution of their node degrees is a linear preferential attachment (LPA). It is applied both to directed and undirected graphs, see Norros and Reittu (2006), Wan et al., (2020) 
among others. The attachment of new nodes starts from an initial network. A newly appended node may connect to $m_0\ge 1$ existing nodes. By the LPA a newly appended node chooses an existing node $i$ randomly in the graph $G=(V,E)$ ($V$ and $E$ denote sets of nodes and edges) with a probability proportional to its degree $k_i$ (i.e. the number of edges of node $i$): $P_{PA}(i)=k_i/\sum_{j\in V}k_j$  (Bagrow \& Brockmann, 2012).
The LPA models provide a "rich-get-richer" mechanism since the earliest appended nodes get likely more edges. This leads to a power-law node degree distribution with tail index $\beta>0$ and a 
constant $C>0$:
\begin{eqnarray*}
 && P\{k_i=j\} \sim  Cj^{-1-\beta},~~ j\to\infty. 
\end{eqnarray*}
Here, 
$a_n\sim b_n$  means that the sequences $a_n$ and $b_n$ are asymptotically equal, i.e., $a_n/b_n \to 1$, $n \to \infty$. The distribution tail becomes heavier for a smaller $\beta$.

However, some real world networks under the constraints of
geographical conditions or topological characteristics have exponential degree distribution, such
as the North American Power Grid Network  (Albert \& Barab\'asi, 2002),
the Worldwide
Marine Transportation Network (Guimer\'{a} et al., 2003)
and some email networks (Deng, 2009).

In the following we focus at the clustering attachment (CA) proposed in  Bagrow and Brockmann (2012)
that constitutes another approach to model the evolving networks. 
 The CA has not received a further development since it seems to be  applicable to  local networks corresponding to a close geolocation of nodes. The CA cannot model the evolution in real-world networks where rare giant nodes may rapidly grow. Instead of giant nodes dense communities are generated around new nodes since nodes are drawn not towards hubs, but towards densely connected groups (Bagrow \& Brockmann, 2012).
 Transport networks may be mentioned as one of practical potential application of the CA. Really, the appearance of a new metro station in a megapolis gives rise to a rapid development of the district and generates a dense transport hub and infrastructure around. Over time, the district reaches a stable state and does not develop so fast anymore. Our paper is devoted to study  the CA including the investigation of the impact of node and edge deletion  which is a novelty.


The CA model  can be described by a graph sequence $G_t$, $t=0,1,\dots$. 
Let us denote the set of nodes and edges in the graph $G_t=(V_t, E_t)$ at step 
$t$, $t \ge 0$ as $V_t$ and   $E_t$, respectively, and their cardinality as $|| V_t||$ and $|| E_t||$.
The evolution starts with a finite undirected initial graph  $G_0$ which consists of  $|| V_0|| \ge 1$ nodes and \textbf{$|| E_0||\ge 0$} edges.  We form $G_{t+1}$ from $G_t$ for $t\ge 0$ by the following rules:

(i)  We append a new node $||V_0||+t+1$ to $G_t$.

(ii) By using a weighted sampling without replacement (see, e.g. Devroye (1986))
we choose a set of nodes  $W_t \subset V_t$, such that $||W_t||=m_0\ge 2$, 
and create  $m_0$ new edges between a newly appended node $||V_0||+t+1$ and $m_0$ selected nodes.  The probabilities 
\begin{equation}\label{CA-norma}
  P_{CA}(i, t) = \frac{c^{\alpha}_{i,t} + \epsilon}{\sum_{j \in V_t} c^{\alpha}_{j,t} + \|V_t\|\epsilon}, \quad i \in V_t, \ \alpha>0
\end{equation}
and 
\begin{equation} \label{CA-norma0}
     P_{CA}(i, t) = \frac{\Ii\{c_{i,t}>0\} + \epsilon}{\sum_{j \in V_t} \Ii\{c_{j,t}>0\} + \|V_t\|\epsilon}, \quad i \in V_t, \ \alpha=0
\end{equation}
are used as the weights for the sampling without replacement.
Here, 
$\Ii\{\cdot\}$ denotes the indicator of the set $\{\cdot\}$ and  $\alpha, \epsilon \geq 0$ 
are attachment parameters.
The clustering coefficient $c_{i,t}$ of node $i \in V_t$ is defined by 
\begin{equation}\label{clustering}
 c_{i,t} = 
 \left\{
  \begin{array}{ll}
    0, & \hbox{$k_{i,t}=0$ or $k_{i,t}=1$,} \\
    2\Delta_{i,t}/\left(k_{i,t}(k_{i,t} - 1)\right), & \hbox{$k_{i,t}\ge 2$.}
  \end{array}
\right.
\end{equation}
Here, $\Delta_{i,t}$ is the number of links between neighbors of node $i$ or, equivalently, the number of triangles involving node $i$, and $k_{i,t}$ denotes the degree of node $i$, both at time $t$. Hereinafter, where it is needed, we will use notation $CA^{(\alpha, \epsilon)}$.

Assumptions on the initial graph $G_0$ depend on $\epsilon$. If $\epsilon>0$ it is enough to assume that
$\|V_0\|\ge m_0$. If $\epsilon=0$, then we have to assume that
$\|\tilde{V}_0\| \ge m_0$,
where the sequence of sets $\tilde{V}_0, \tilde{V}_1, \dots$
is defined as follows: $\tilde{V}_n=\{ i \in V_n: \ c_{i,n}>0\}$,
$n \ge 1$. Let ${\tilde E}_t$ denote the set of all possible edges between  pairs of nodes $\{i_1,i_2\}$, such that $i_1\in {\tilde V}_t$, $i_2\in {\tilde V}_t$, $i_1 \not =i_2$, including probably not existing edges. It is worth noting that  $(\tilde{V}_n, \tilde{E}_n)$ is a complete graph. 


A probability mass function $ P_{CA}(i, t)$, $i\in V_t$ for the  limit case $\alpha\downarrow 0$
is constructed by us keeping in mind the relation $x^{\alpha} \to 1$ as $\alpha \downarrow 0$ for any fixed positive $x$ and by the convention  $0^0=0$.

The CA model  with probabilities
\begin{equation}\label{CA}
  P_{CA}(i,t) \propto c^{\alpha}_{i,t} + \epsilon
\end{equation}
was used in Bagrow and Brockmann (2012).
We recall that  $x\propto y$  means that there is a non-zero constant $C$ such that $x = C \cdot y$. 
In fact, (\ref{CA-norma}) and (\ref{CA}) define the same conditional probabilities 
with $C=M_t$, 
where 
\begin{eqnarray*}
M_t&=&\left(\sum_{j \in V_t} c^{\alpha}_{j,t} + \|V_t\|\epsilon\right)^{-1}.
\end{eqnarray*}
Such kind of the CA excludes the appearance of multiple edges. Furthermore,  the sets $V_t$ and $E_t$ grow linearly in $t$: $||V_t||=||V_0||+t$ and $||E_t||=||E_0||+m_0t$, $t \ge 0$
without the node and edge deletion.
Whence it immediately follows that the ratio $||E_t||/||V_t||$ tends to $m_0$ as $t \to \infty$.

In the paper we consider several  modifications of the  CA. 
One of them allows  a node deletion when a new node is appended. Specifically, together with (i), (ii), 
one more rule is added:

(iii) We choose one node  from the set $V_{t}$ uniformly, i.e. each node with probability $1/||V_{t}||$, and delete it. We remove also all edges belonging to this node.

The second modification of the CA model concerns the allowance of an edge deletion when a new node is appended. Then the rule (iii) 
is replaced by the following:

(iv) We choose one edge from the set $E_{t}$ uniformly, i.e. each edge with probability $1/||E_{t}||$, and delete it. 

Further, we keep the notation CA for   the clustering attachment model generated by rules (i), (ii).
 \par
In contrast to the PA, 
the intuition outlined in  Bagrow and Brockmann (2012) shows 
that the CA  does not 
lead to  a power-law node degree distribution, but to an  exponential light- tailed distribution.
Our {\bf {first objective}} is 
to find lower and upper bounds for the increment of consecutive mean clustering  coefficients over time; to derive that (total) node  triangle counts 
and 
node degrees  are submartingales irrespective of
the CA parameters $\alpha$ and $\epsilon$ in (\ref{CA-norma}), (\ref{CA-norma0}). 
In Markovich and Vai{\v c}iulis (2024)
a special attention was devoted to the model $CA^{(0,0)}$.  It  was proved there without additional assumptions that a total triangle count tends to infinity almost surely for the latter model.
\par
Our {\bf {second objective}} is  to check the hypothesis of a light-tailed node degree distribution with regard to 
$\alpha$, $\epsilon$ in (\ref{CA-norma}) 
as well as to  evolution strategies, namely, without node and edge deletion and with the deletion of one of the existing nodes or edges each time when a new node is appended. 
This is done by the  extreme value index  (EVI) estimation of the node degree, see Section \ref{2.4} for the EVI definition. The EVI allows us to distinguish between light- and heavy-tailed distributions. For heavy-tailed distributions, the value of the EVI indicates the heaviness of the distribution tail. The EVI  will be estimated by several known semi-parametric estimators. Moreover, our simulation study aims to confirm the results obtained in the propositions.

The paper is organized as follows. Related definitions and results are given in Sect. \ref{Sec2}. Theoretical results are provided in Sect. \ref{Sec3}. The simulation study is presented in Sect. \ref{Sec4}. Sect. \ref{Sec5} illustrates our theoretical results for real networks.
We finalize with the conclusions in Sect. \ref{Sec6}. 

%% file: Sec2ANOR.tex
\section{Related definitions and results}\label{Sec2}
\subsection{Mean clustering coefficient
}\label{sub-CA}

The clustering coefficient $c_{i,t}$ of node $i \in V_t$ in (\ref{clustering})  measures its tendency  to form triangles of the nearest nodes in its neighborhood. The average clustering coefficient for an undirected graph $G_t = (V_t, E_t)$, $t \ge 0$ is defined  by
\begin{equation}\label{av-clust}
  \overline{C}_t = \frac{1}{||V_t||}\sum_{i \in V_t}c_{i,t}.
\end{equation}
For the CA model the definition (\ref{av-clust}) can be rewritten as follows:
$$
  \overline{C}_t = \frac{1}{||V_0||+t}\sum_{i=1}^{||V_0||+t} c_{i,t}.
$$
The statistic $\overline{C}_t$ may be used to  distinguish
between geometric and non-geometric networks.  We refer to Bringmann et al. (2019); Michielan et al. (2022)
for the explanation of the last two notions.  Roughly speaking, in geometric networks each node is positioned in a geometric space and there is a distance to other nodes that can be measured.
Typically, if the value
of $\overline{C}_t$ does not vanish in $t$, then this is considered to be evidence for a geometry.
For example,  the  average clustering coefficient
of a geometric inhomogeneous random graph  does not vanish as $t$ increases (Bringmann et al., 2019).
For  hyperbolic
random graphs,  $\overline{C}_t$ converges in probability to a positive constant (Fountoulakis et al., 2021).
For the inhomogeneous random graphs the decay of  $\overline{C}_t$ can be extremely slow 
(Michielan et al., 2022). 

\subsection{Extreme value index}\label{2.4}
Let us recall the definition of the EVI (de Haan \&  Ferreira, 2007).
A quantile type function $U$ associated with an arbitrary  cumulative distribution function (cdf)
$F(x)$  is defined by
$$
U(t)= \left\{
  \begin{array}{ll}
    0, & \hbox{$0<t\le 1$,} \\
    \inf\left\{x: \ F(x) \ge 1-(1/t)\right\}, & \hbox{$t>1$.}
  \end{array}
\right.
$$
The function $U(t)$ is said
to be of an extended regular variation if there exists a positive function $a(t)$ such that for some $\gamma \in R$ 
and
all $x>0$, it holds
$$
     \lim_{t \to \infty} \frac{U(tx)-U(t)}{a(t)}=
    \left\{
  \begin{array}{ll}
    \ln(x), & \hbox{${\gamma}=0$,} \\
    (x^{\gamma}-1)/{\gamma}, & \hbox{${\gamma} \not =0$.}
  \end{array}
\right.
$$
The parameter $\gamma$ is called the EVI. If $\gamma>0$ holds,
then the EVI and tail index $\beta$ are related by $\beta=1/\gamma$. For example, EVI of Burr or Fr\'{e}chet distribution is positive.
The examples of light-tailed distributions with $\gamma=0$
are an exponential distribution and  a normal distribution, while the examples of light-tailed distributions with  $\gamma<0$ are uniform and beta distributions. 
It is worth mentioning  that
not all distributions belong to the class of extended regular variation.  
A Poisson distribution can be taken as an example. 

%% file: Sec3ANOR.tex
\section{Main Results}\label{Sec3}
\par

In this section we provide our theoretical results related to the CA model. The model is considered without node and edge deletion here.

\input{subsec3.1ANOR}
\input{subsec3.2ANOR}
\input{subsec3.3ANOR}

%% file: subsec3.1ANOR.tex
\subsection{Weighted sampling without replacement}\label{Sec3.1}

Let us adopt  the sampling procedure for our purposes. Suppose that a 
random graph $G_t=(V_t, E_t)$ is generated after $t$ evolution steps. We
denote the unordered pair  of nodes chosen from $V_t$ by using
weighted sampling without replacement (WSwR)  as  $W_t$. Let
a set $\mathcal{E}_t$ be such, that $(V_t, \mathcal{E}_t)$ forms a complete graph.
The attachment probabilities
\begin{eqnarray*} 
&& P\left(W_t=\{i_1, i_2\}\right), \quad \{i_1, i_2\} \in \mathcal{E}_t
\end{eqnarray*}
are needed to generate the $CA^{(\alpha, \epsilon)}$ model.

\begin{proposition}\label{prop1} Let the sequence of random graphs $G_1,G_2, \dots$ be generated by the $CA^{(\alpha,
\epsilon)}$ model with the initial graph $G_0$ and  parameters $\alpha> 0$, $\epsilon\ge 0$ and $m_0=2$.
Then for any $\{i, j\} \in \mathcal{E}_t$,
\begin{equation} \label{w01}
P(W_t=\{i, j\}) = \frac{P_{CA}(i,t) P_{CA}(j,t) \left(2-P_{CA}(i,t)
-P_{CA}(j,t)
\right)}{\left(1-P_{CA}(i,t)\right)\left(1-P_{CA}(j,t)\right)}.
\end{equation}
\end{proposition}

\begin{corollary} \label{cor-01} 
Depending on  $\epsilon>0$ or $\epsilon=0$  the collection of conditional probabilities (\ref{w01}) form a conditional probability distribution on $\mathcal{E}_t$ and
${\tilde E}_t$, respectively.
\end{corollary}

%% file: subsec3.2ANOR.tex
\subsection{
Bounds for increments of the mean clustering coefficients
}


We investigate
the following sequences of a fixed node
$i\in V_t$.
Namely, we consider node
degrees  $\{k_{i, t}\}$,
the triangle counts $\{\Delta_{i, t}\}$ and the clustering coefficients $\{c_{i, t}\}$ for all $t\ge 0$.
%
For any $t \ge 1$ we assume 
\begin{equation} \label{df00}
    k_{||V_0||+t,s}=0, \quad 0 \le s <t.
\end{equation}
The same assumption holds for the rest of the sequences.
Note that the sequence $\{k_{i, t}\}_{t \ge 0}$ is non-decreasing, i.e., $k_{i,t} \le k_{i,t+1}$ for any $t \ge 0$. The 
same is valid for  $\{\Delta_{i, t}\}_{t \ge 0}$. As for $\{c_{i, t}\}_{t \ge 0}$, it is stated in
Bagrow and Brockmann (2012) 
that "even when a new triangle is formed, the clustering coefficient after an 
attachment is almost always less than it was before". Here, the notion "almost always" is used in the sense that it might be $c_{i,t} > c_{i,t+1}$
for some $t \ge 0$.
 Let us consider the following counter examples.
 \begin{example}\label{Exam3.1a}
 Let $V_0=\{1,2,3,4\}$ hold
and the graph $G_0$ be a rectangle.
Let the numbering of nodes go clockwise.
Assuming $m_0=2$,  a newly appended node is connected with nodes $1$ and $3$ at each time step $t \to t+1$, $t\ge 0$.
Under such evolution we have $\Delta_{i,t}=0$ and consequently, $c_{i,t}=0$ for any $i \in V_t$ and $t \ge 0$.
\end{example}
The next simple modification of the initial graph $G_0$ demonstrates that the clustering coefficients of some nodes can remain constant.
\begin{example}\label{Exam3.2a}
Let $G_0$ be the same rectangle as in the previous example but with diagonal connecting nodes $1$ and $3$. Under the same evolution as in the previous example, we have
$c_{2,t}=c_{4,t}=1$ for all $t \ge 0$ by (\ref{clustering}).
\end{example}
The last two examples demonstrate the impact of the initial graph $G_0$ on the average clustering coefficient ${\bar C}_t$, see  (\ref{av-clust}) for its definition. Turning back to Example 
\ref{Exam3.1a}
we have ${\bar C}_t=0$ for any $t \ge 0$.
In Example \ref{Exam3.2a} we have
 by 
 $$
c_{1,t}=c_{3,t}=\frac{2}{t+3}, \quad  c_{2,t}=c_{4,t}=1
$$
and $c_{i,t}=1$ for $t \ge 1$, $||V_0||<i \le ||V_0||+t$, where $||V_0||=4$, that
${\bar C}_t=4/((t+3)(t+4))+(t+2)/(t+4)$ holds.
Whence it follows ${\bar C}_t \to 1$ as $t \to \infty$.

Let us consider the behaviour of ${\bar C}_t-{\bar C}_{t+1}$ in a general case.

\begin{proposition}\label{Prop1}
Let $G_t=(V_t, E_t)$, $t \ge 0$ be the sequence of the CA evolved graphs.
Assume that the evolution parameter $m_0$ satisfies $m_0=2$.
Then for any $t \ge 0$, it holds
\begin{equation} \label{ineqa-01}
-\frac{3}{||V_0||+t+1} \le {\bar C}_t -{\bar C}_{t+1}  \le
\frac{7/3}{||V_0||+t+1}. 
\end{equation}
\end{proposition}
We have several notices related to Prop. \ref{Prop1}.
\begin{remark}
\begin{enumerate}
\item Proposition \ref{Prop1} can be generalized to the case $m_0 \ge 2$.

\item
The inequality 
(\ref{ineqa-01}) is 
true for any realization of the $CA^{(\alpha,\varepsilon)}$
model. Since we do not use distribution (\ref{CA-norma}) in our proof, (\ref{ineqa-01}) 
holds true if we replace the CA model by the LPA model.

\item

Let $0\le t <t'$ hold. By using the inequality
$|{\bar C}_{t'} -{\bar C}_{t}| \le \sum_{s=t}^{t'-1} \left| {\bar C}_{s+1} -{\bar C}_{s} \right|$
and  (\ref{ineqa-01}) it follows that
$$|{\bar C}_{t'} -{\bar C}_{t}| \le  \sum_{s=t}^{t'-1} \frac{3}{||V_0||+s+1}\le
3 \left(H(t')-H(t) \right),$$
where $H(t)=\sum_{s=1}^{t-1} 1/s$.

By using the inequalities $\ln(t-1) \le H(t) \le 1+\ln(t-1)$, $t\ge 2$ we get
$$
|{\bar C}_{t'} -{\bar C}_{t}| \le 3\ln\left( \frac{{\rm e} (t'-1)}{t-1} \right),
$$

where ${\rm e} \approx 2.718$ is Euler's number.
The last inequality does  not imply that
the limit
of the sequence ${\bar C}_t$, $t \ge 0$ exists.
\item From (\ref{ineqa-01}) it  follows immediately that
a sequence 
$\{{\bar C}_t-{\bar C}_{t+1}\}_{t \ge 0}$ converges to zero as $t\to\infty$.
\end{enumerate}
\end{remark}

\begin{remark}\label{remark-monotone} The proof of Proposition \ref{Prop1} implies that ${\bar C}_{t}$, $t \ge 0$  is not a monotone sequence.
Monte Carlo simulation shows
that 
${\bar C}_{t}$ tends to zero as $t\to\infty$, see Fig. \ref{fig:4.8.2}.
\end{remark}

We cannot extend Prop. \ref{Prop1} to graphs
such that a node or an edge can be deleted at each step of the CA evolution  in a similar way. The reason is the following. If   an edge   between some nodes $i,i'\in V_t$ (or one of the nodes) 
is deleted, one has to take into account how this deletion
affects on the rest of the graph.

%% file: subsec3.3ANOR.tex
\subsection{Submartingality} \label{nd}

We introduce a sequence of $\sigma$-algebra's as follows:
\begin{eqnarray*}&&\mathcal{F}_t=\sigma\left\{G_0, G_1, \dots, G_t\right\}, \quad t=0,1,2,\dots
\end{eqnarray*}
Obviously, it holds $\mathcal{F}_t \subset \mathcal{F}_{t'}$ when $t \le t'$, i.e.,
the sequence  $\mathcal{F}_0, \mathcal{F}_1, \dots$ forms a filtration.

\begin{proposition}\label{Prop2} Let the assumptions of Prop. \ref{prop1} 
be satisfied. Then for any fixed node $i$,
the sequence $(k_{i,0},\mathcal{F}_0), (k_{i,1},\mathcal{F}_1), \dots$ is  a submartingale.
\end{proposition}

\begin{proposition}\label{Prop3} Let the assumptions of Prop. \ref{prop1} be satisfied. Then for any fixed node $i$,
the sequence $(\Delta_{i,0},\mathcal{F}_0), (\Delta_{i,1},\mathcal{F}_1), \dots$ is  a submartingale.
\end{proposition}

Let us recall that the total count of triangles in the graph $G_n$ is defined by
\begin{equation}\label{df04}
\Delta_t=\frac{1}{3} \sum_{i \in V_t} \Delta_{i,t}.
\end{equation}

\begin{corollary}\label{Cor3.2}
Let the assumptions of Prop. \ref{prop1} be satisfied. Then
the sequence 
\begin{eqnarray*}&&(\Delta_{0},\mathcal{F}_0), (\Delta_{1},\mathcal{F}_1), \dots\end{eqnarray*} is  a submartingale.
\end{corollary}

%% file: Sec4ANOR.tex
\section{Simulation study}\label{Sec4}
\par
Our simulation study aims to support the theoretical results obtained in Prop. \ref{Prop1}-\ref{Prop3} for the $CA^{\alpha,\epsilon}$ model without node and edge deletion and to investigate the CA evolution with a uniform node or edge deletion.
\par
Let $G_t=(V_t, E_t)$ be the graph at evolution step $t$.  We  consider the CA evolution with  the creation of
$m_0\ge 2$
new edges  at each evolution step when a new node is appended. 

\subsection{Illustration of the CA evolution with and without node and edge deletion}
\par
\begin{figure}[tbp]
 \begin{minipage}[t]{\textwidth}
\centering
  \subfigure[]{\includegraphics[width=0.22\textwidth]{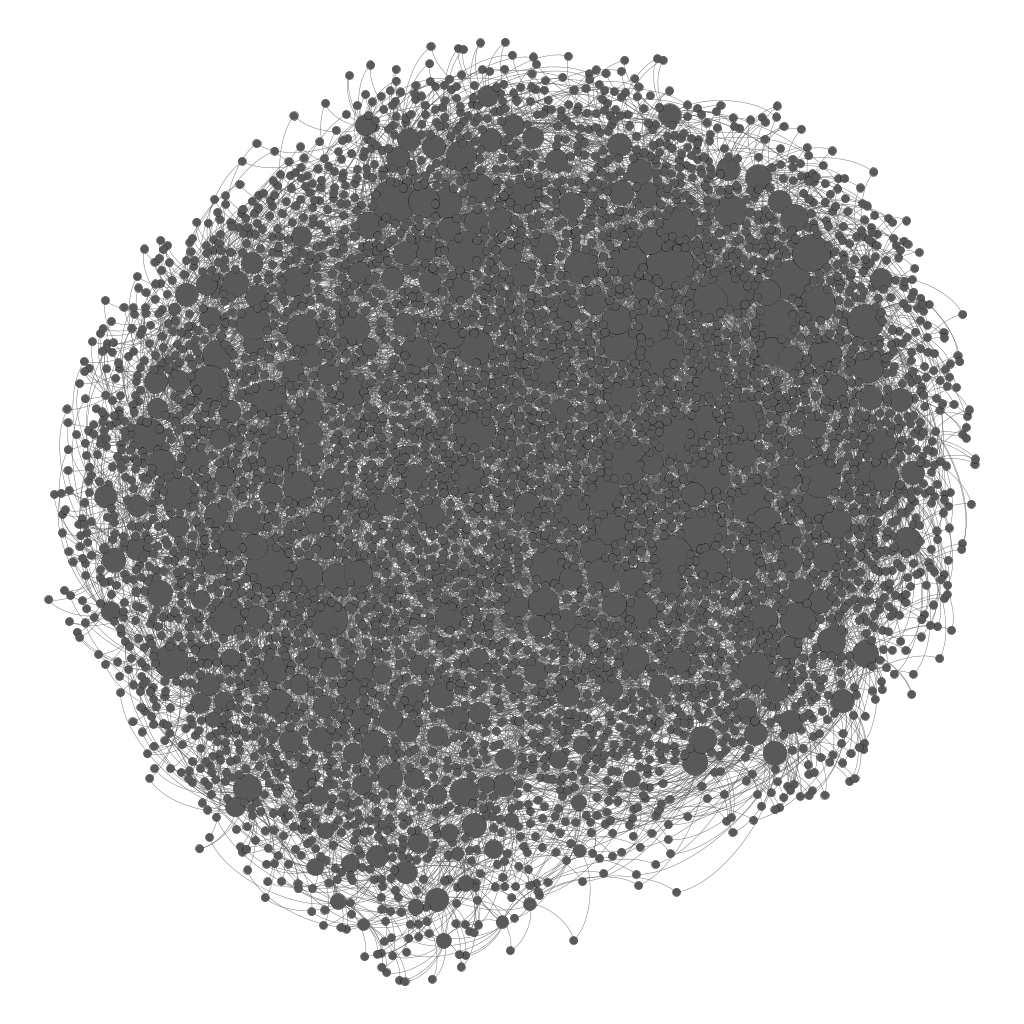}
                 \label{fig:1a}}
  \subfigure[]{\includegraphics[width=0.22\textwidth]{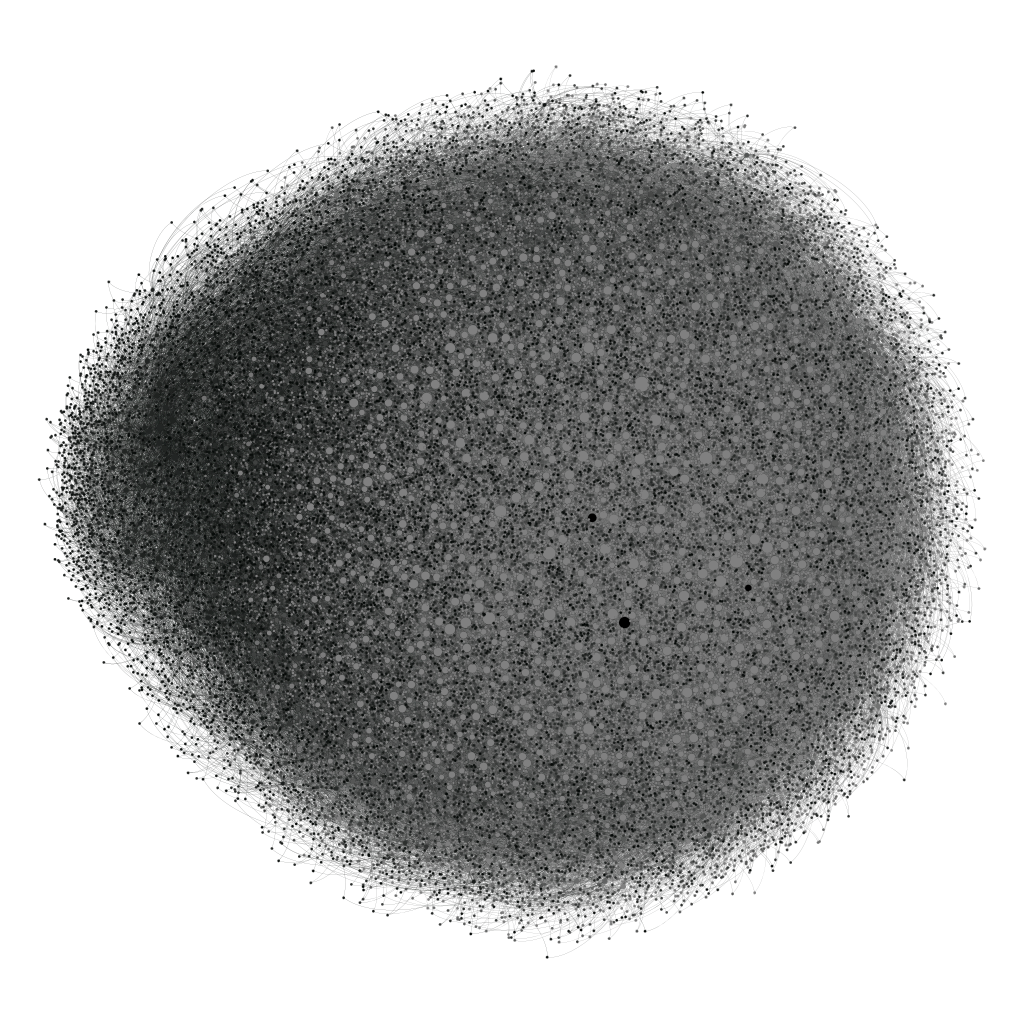}
          \label{fig:1b}}
 \subfigure[]{\includegraphics[width=0.22\textwidth]{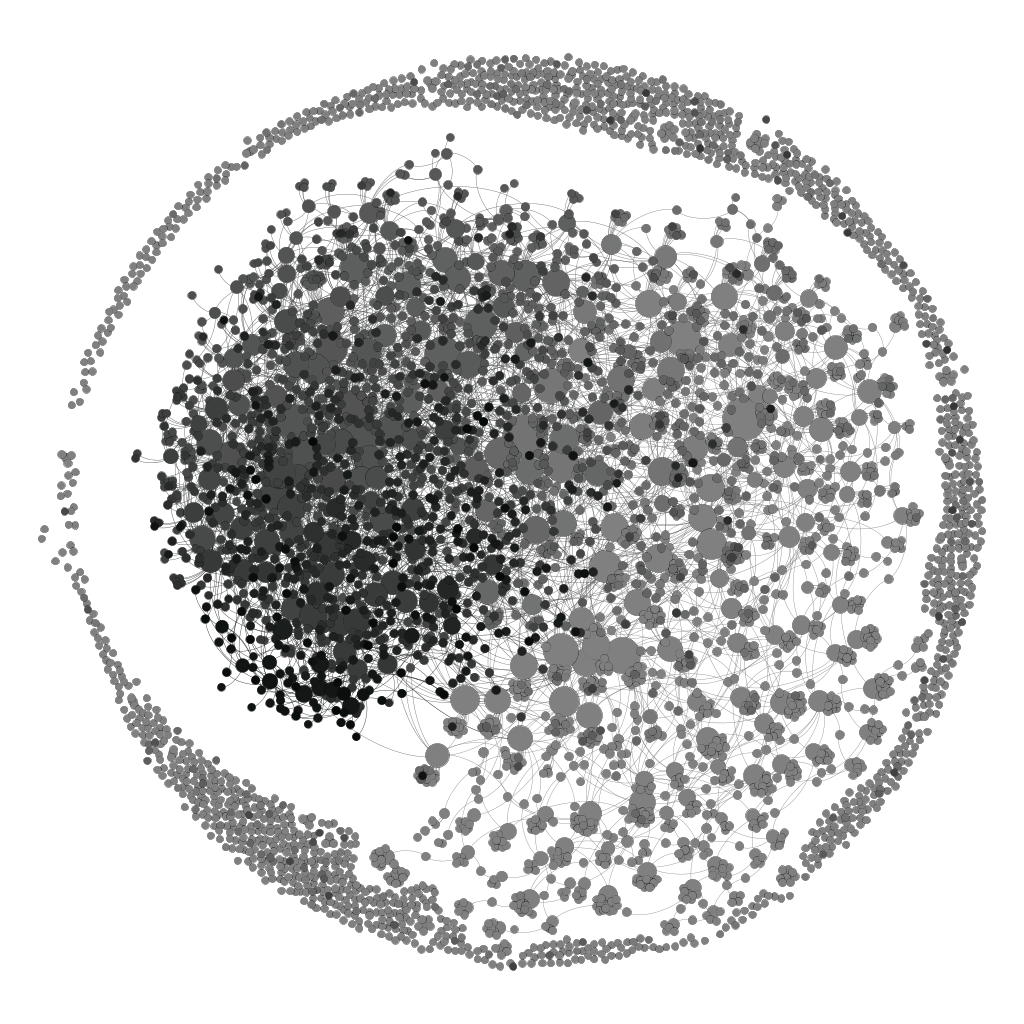}
          \label{fig:1c}}
\subfigure[]{\includegraphics[width=0.22\textwidth]{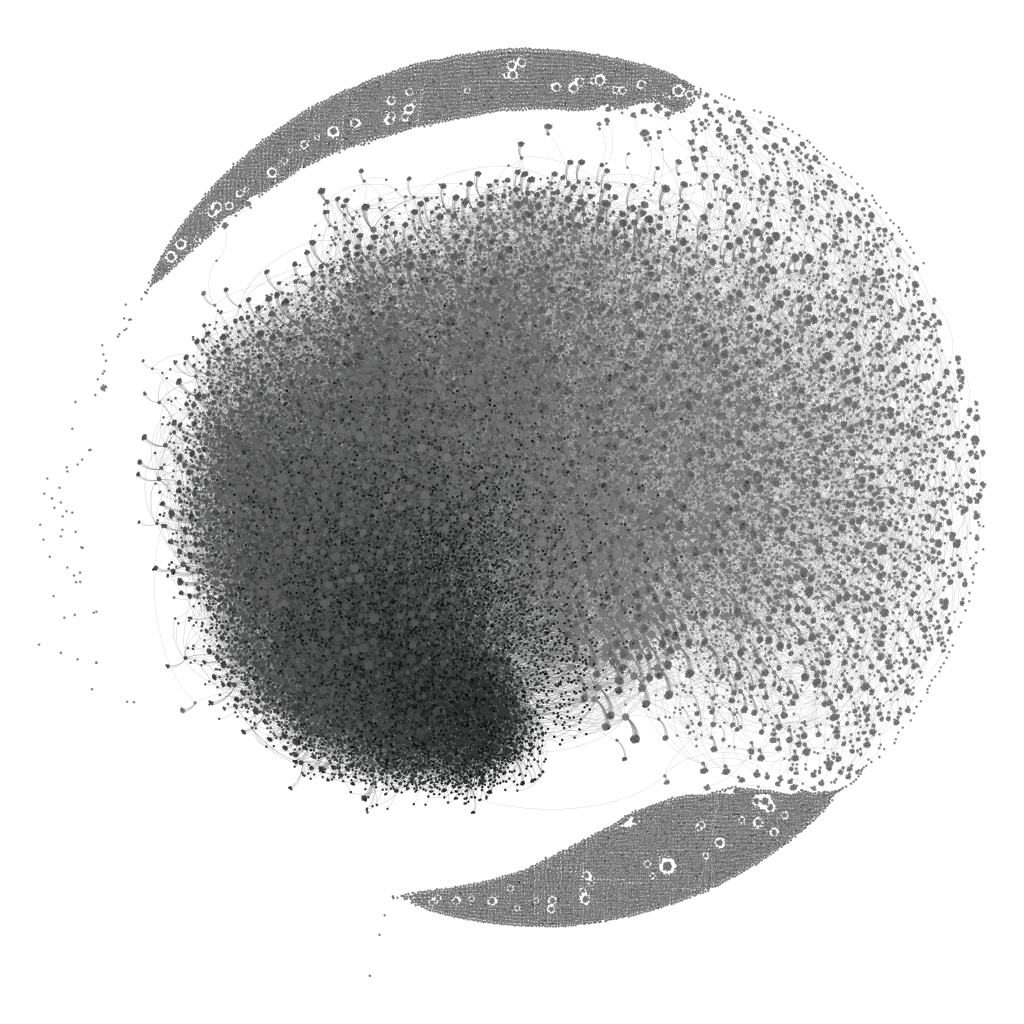}
          \label{fig:1d}}
        \caption{The graph obtained by the CA with attachment probability  (\ref{CA-norma}) and  parameters $(\alpha, \epsilon) =(1,0)$ 
        and $m_0=2$  from the initial graph in Fig. \ref{fig:1a} after $5\cdot10^4$ evolution steps without node and edge deletion (Fig. \ref{fig:1b}), with a uniform node deletion (Fig. \ref{fig:1c}) and with a uniform edge deletion (Fig. \ref{fig:1d}). The node size is proportional to the node degree, and the node color represents the "life time" of the node,
    i.e. the "older" the node the darker the color.
       }\label{fig:1} 
       \end{minipage}
       \end{figure}

An example of an initial graph $G_0$
obtained by the CA  with parameters $m_0=2$, $(\alpha,\epsilon)=(1,0)$ starting from a triangle of connected nodes and containing $5\cdot 10^3$ nodes is shown in Fig. \ref{fig:1a}.
The evolution without node and edge deletion leads to   nodes with a dense  structure, see  Fig. \ref{fig:1b}. A uniform node or  edge deletion  at each evolution step when a new node and new edges are appended,  causes the appearance of isolated nodes. The proportion of the latter nodes is $31.4\%$  and $26.3\%$  for uniform node  and  edge deletion in Fig. \ref{fig:1c}, \ref{fig:1d}, respectively. Naturally, 'the removal of a node implies the malfunctioning of all its edges as well, node removal inflicts more damage than edge removal' (Albert \& Barab\'asi, 2002).
If
nodes are isolated, then the attachment to them is unlikely since their triangle counts are zero-valued.

       \subsection{Averaging of total triangle counts}
 \par
       \begin{figure}
    \begin{minipage}[t]{\textwidth}
\centering
  \subfigure[]{\includegraphics[width=0.31\textwidth]{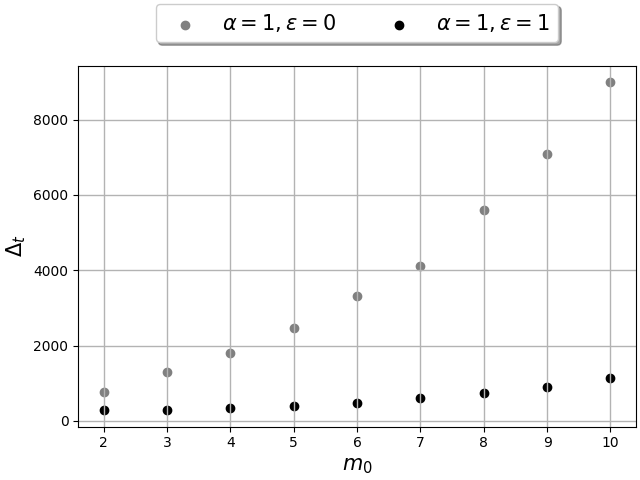}
                  \label{fig:11a}}
    \subfigure[]{\includegraphics[width=0.31\textwidth]{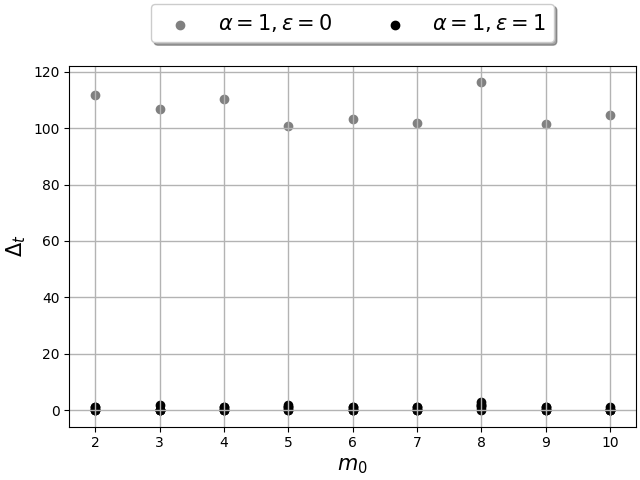}
                   \label{fig:11b}}
    \subfigure[]{\includegraphics[width=0.31\textwidth]{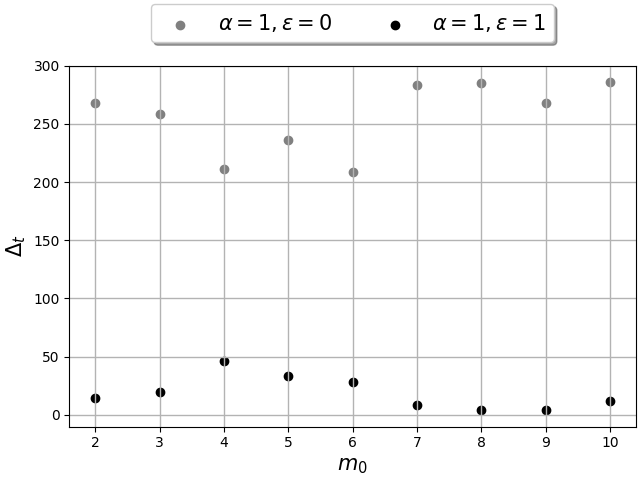}
                    \label{fig:11c}}
    \caption{Averages of total triangle counts $\Delta_t$ over $10$ graphs against $m_0$ for the  CA graphs after $t=5\cdot10^4$ evolution steps with sets of parameters $(\alpha, \epsilon) \in \{(1, 0), (1, 1)\}$ without node and edge deletion (Fig. \ref{fig:11a}), with a uniform node deletion (Fig. \ref{fig:11b}) and with a uniform edge deletion (Fig. \ref{fig:11c}).}
    \label{fig:11}
    \end{minipage}
\end{figure}
 Fig. \ref{fig:11a} shows that the mean triangle counts grow with increasing $m_0$ for the CA model  without node and edge deletion for both pairs of parameters  $(\alpha,\epsilon) \in \{(1,0), (1,1)\}$.
 The case $(1,0)$ provides a higher increasing rate  than $(1,1)$ since  the impact of the clustering coefficient for $(1,1)$ is by (\ref{CA-norma}) weaker than for $(1,0)$ and   node $i$ can likely get a new edge uniformly irrespective of the number of triangles involving this node.
In contrast, the evolution with node or  edge deletion in  Fig. \ref{fig:11b}, \ref{fig:11c} leads to a fluctuation around  constants for the average triangle counts. For $(\alpha,\epsilon) = (1,1)$ the triangles  appear rarely or they are absent in the graph for any considered values of $m_0$. For $(\alpha,\epsilon) = (1,0)$ the triangle counts are approximately constant with regard to $m_0$, and their values are larger for the evolution with edge deletion rather than for the evolution with node deletion.

\subsection{The extreme value index estimation}

\begin{figure}[ht!]
\begin{minipage}[t]{\textwidth}
         \centering
   \subfigure[]{\includegraphics[width=0.45\textwidth]{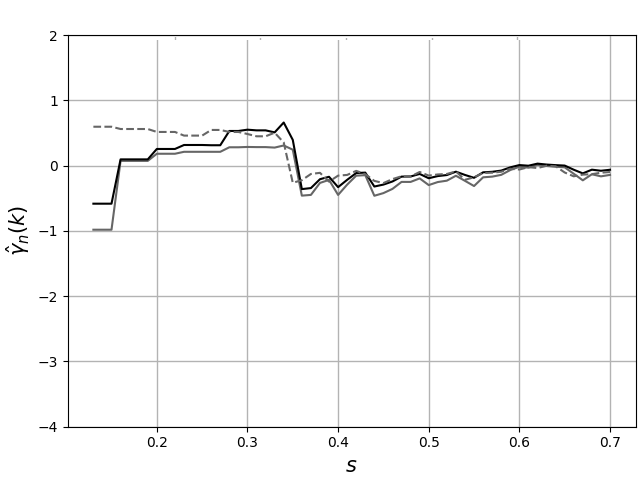}
    \includegraphics[width=0.45\textwidth]{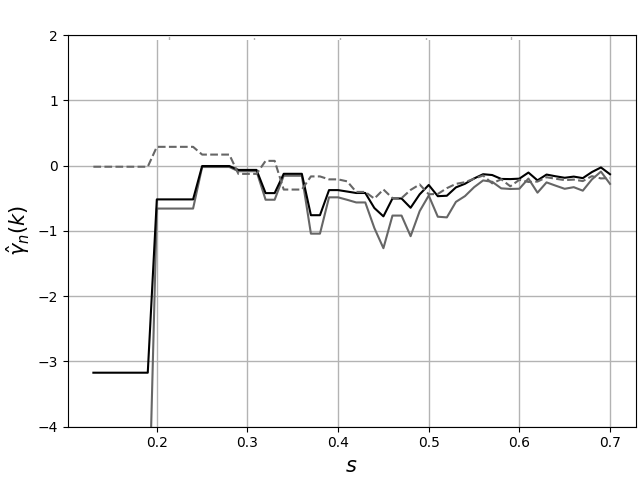}
                \label{fig:7a}}
    \subfigure[]{\includegraphics[width=0.45\textwidth]{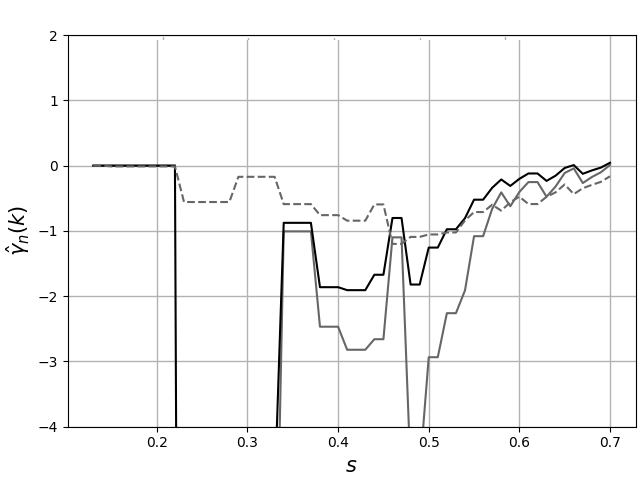}
   \includegraphics[width=0.45\textwidth]{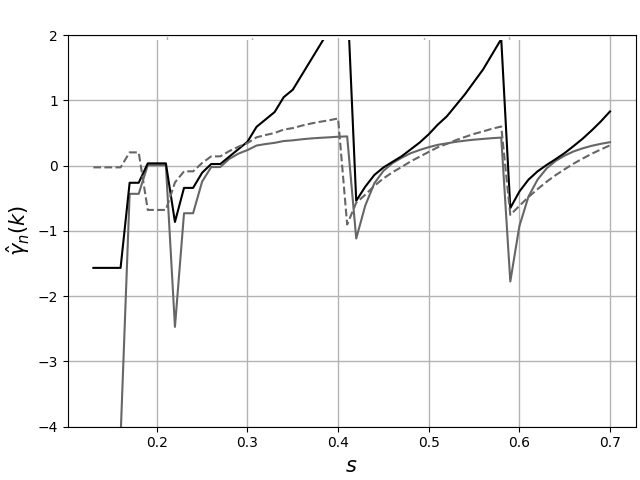}
             \label{fig:7b}}
    \subfigure[]{\includegraphics[width=0.45\textwidth]{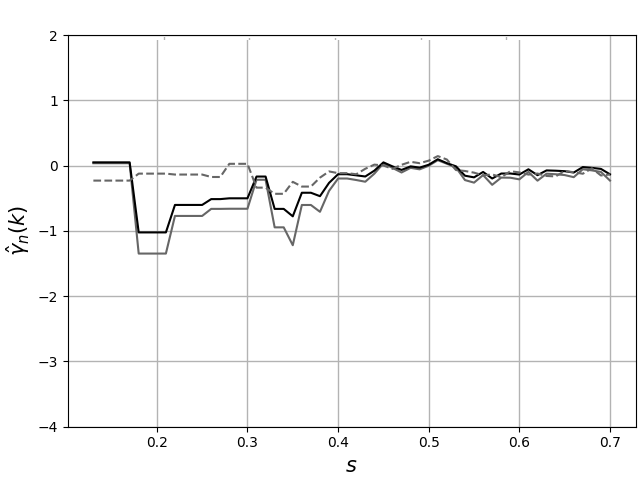}
     \includegraphics[width=0.45\textwidth]{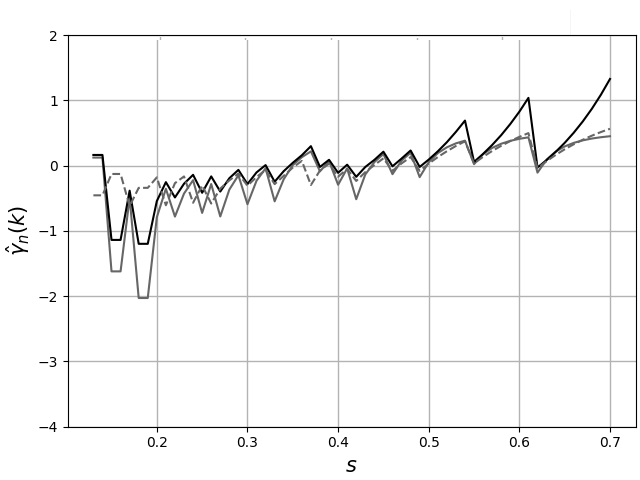}
                 \label{fig:7c}}
       \caption{EVI plots of the moment (grey line), mixed moment (black line) and UH (dotted line) estimates $\hat\gamma_n(k)$ of node degrees against the parameter $s$ included in the number of the largest order statistics $k=[n^{s}]$, $n = \|V_t\|$, for graphs evolved by the CA with probability (\ref{CA-norma}) and parameters $(\alpha, \epsilon)=(1,0)$ (the left column), 
    and $(\alpha, \epsilon)=(1,1)$ (the right column) and
    with $m_0=2$ after $t=5\cdot10^4$ evolution steps  without node and edge deletion
    (Fig. \ref{fig:7a}), with a uniform node deletion (Fig. \ref{fig:7b}) and with a uniform edge deletion
    (Fig. \ref{fig:7c}).
        }\label{fig:7}
        \end{minipage}
\end{figure}
\begin{figure}[ht!]
\begin{minipage}[t]{\textwidth}
   \centering
  \subfigure[]{\includegraphics[width=0.45\textwidth]{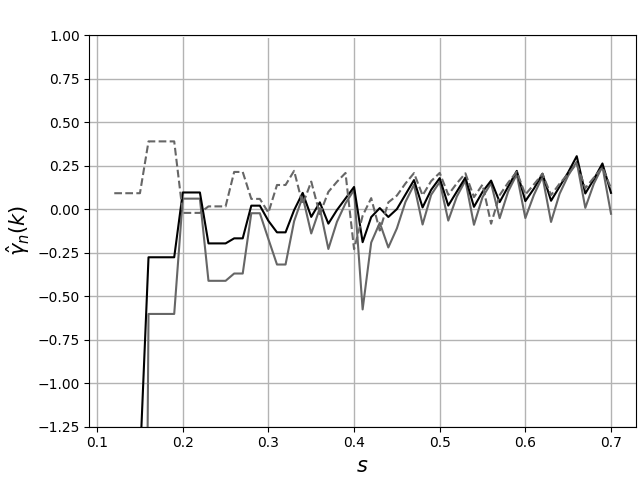}
            \includegraphics[width=0.45\textwidth]{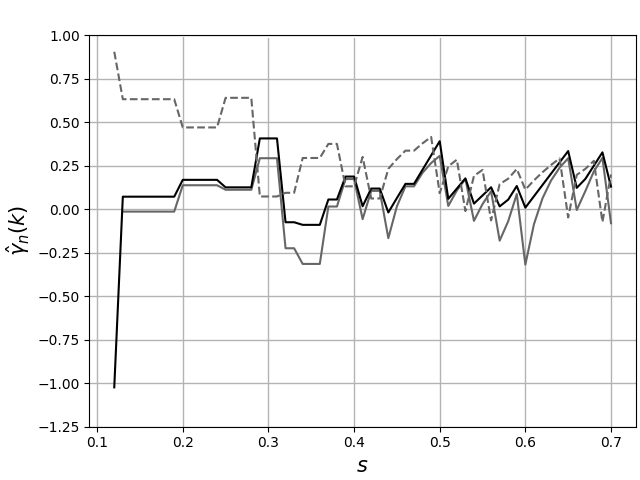}
                 \label{fig:tr-a}}
    \subfigure[]{\includegraphics[width=0.45\textwidth]{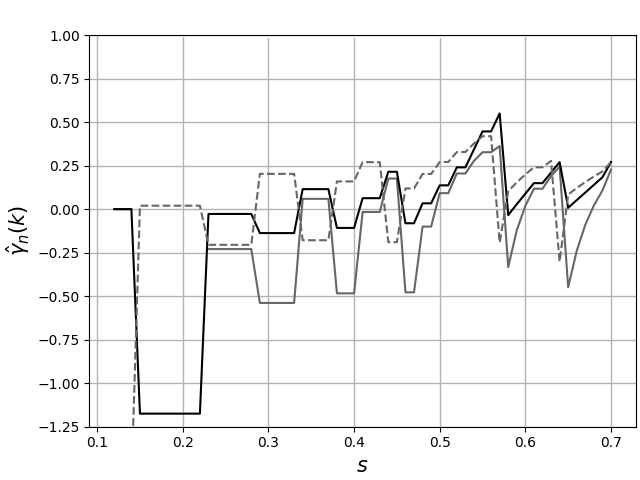}
                 \label{fig:tr-b}}
    \subfigure[]{\includegraphics[width=0.45\textwidth]{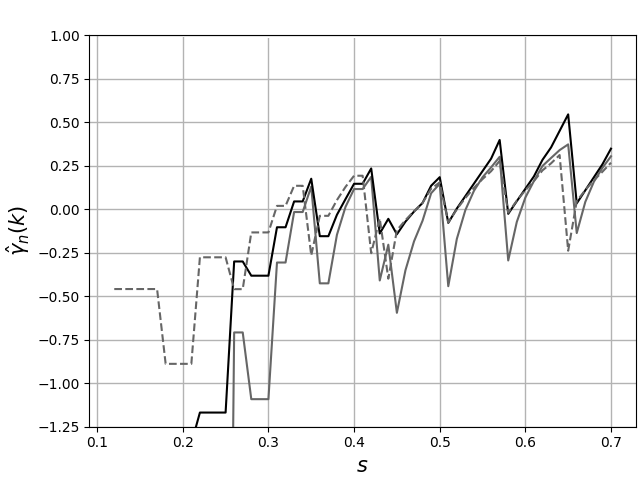}
                  \label{fig:tr-c}}
       \caption{EVI plots of the moment (grey line), mixed moment (black line) and UH (dotted line) estimates $\hat\gamma_n(k)$
    of node triangle counts against the parameter $s$ included in the number of the largest order statistics $k=[n^{s}]$, $n = \|V_t\|$, for graphs evolved by the CA with probability (\ref{CA-norma}) and $m_0=2$ after $t=5\cdot10^4$ evolution steps  without node and edge deletion
    (Fig. \ref{fig:tr-a}), with a uniform node deletion
    (Fig. \ref{fig:tr-b}) and with a uniform edge deletion
    (Fig. \ref{fig:tr-c}). Parameters $(\alpha, \epsilon)$ are equal to $(1,0)$ (left) and to $(1,1)$ (right) in Fig. \ref{fig:tr-a} , and to $(1,0)$ both in Fig. \ref{fig:tr-b}, \ref{fig:tr-c}.
        }\label{fig:tr}
        \end{minipage}
\end{figure}
\begin{table}[h!]
\centering
\addtolength{\tabcolsep}{-5pt}
\begin{tabular}{||c | c c c| c c c| c c c||}
 \hline
 & \multicolumn{9}{c|}{The  evolution} \\
 \hline
 $(\alpha, \epsilon)$ & \multicolumn{3}{c|}{ without node/edge deletion} & \multicolumn{3}{c|}{ with node deletion} & \multicolumn{3}{c||}{ with edge deletion} \\
  & $\|V_t\|$ & $\{k_i > 2\}$ & $\{\Delta_{i,t} > 0\}$ & $\|V_t\|$ & $\{k_i > 0\}$ & $\{\Delta_{i,t} > 0\}$ &  $\|V_t\|$ & $\{k_i > 0\}$ & $\{\Delta_{i,t} > 0\}$\\ 
 \hline
 (1, 0)  & 55000 & 1277 & 1277 & 5000 & 125 & 125 & 55000 & 624 & 624 \\
 (1, 1)    &  55000    & 36757 & 311 & 5000 &  4596 &  2 & 55000 & 5469 & 0 \\
 \hline
\end{tabular}
\caption{The number of nodes $\|V_t\|$ after $t=5\cdot10^4$ evolution steps, where the initial graph in Fig. \ref{fig:1a} contains $5\cdot10^3$ nodes; the number of nodes with degrees larger than the minimum value, i.e. $\{k_{i,t} > k_{min}\}$, $k_{min} \in \{0, 2\}$; and the number of nodes involved in triangles for the CA graphs, i.e. such that $\{\Delta_{i,t} > 0\}$ for  parameters  $(\alpha,\epsilon) \in \{(1,0), (1,1)\}$, $m_0=2$.
}
\label{table:1}
\end{table}
It is  claimed in Bagrow and Brockmann (2012)
that networks evolving by the CA exhibit an exponential tail of the node degree distribution, i.e. the EVI of
the node degree distribution is equal to zero. In order to verify this claim we apply the moment, the mixed moment and the UH  estimators, see (Dekkers et al., 1989; Fraga Alves et al., 2009; Beirlant et al., 2004),
respectively. The latter semiparametric estimators are designed to estimate real-valued EVI.
They are defined by means of the Hill's estimator
\begin{eqnarray*}
&&\hat{\gamma}^H(n,k)=\frac{1}{k}\sum_{i=1}^k 
\ln\left( \frac{X_{(n-i+1)}}{X_{(n-k)}}\right) \end{eqnarray*}
built by the order statistics $X_{(1)}\le X_{(2)}\le ...\le X_{(n)}$ of the observations and for some $k=1,...,n-1$ (Hill, 1975).
The moment estimator is defined as
$$\hat{\gamma}^M_{n,k}=\hat{\gamma}^H(n,k)+1-\frac{1}{2}\left(1-(\hat{\gamma}^H(n,k))^2/S_{n,k})\right)^{-1},
$$ where $S_{n,k}=(1/k)\sum_{i=1}^k\left(\log X_{(n-i+1)}-\log
X_{(n-k)}\right)^2$.
The UH estimator is
\begin{eqnarray*}&&\hat{\gamma}^{UH}_{n,k} =(1/k)\sum_{i=1}^k\log UH_i-\log
UH_{k+1},\end{eqnarray*}
where $UH_i=X_{(n-i)}\hat{\gamma}^H(n,i)$.
The mixed moment estimator is defined as
\begin{eqnarray*}
&&\widehat{\gamma}_{n}^{M M}(k) =\frac{\widehat{\varphi}_{n}(k)-1}{1+2 \min \left(\widehat{\varphi}_{n}(k)-1,0\right)},
\end{eqnarray*}
where
$$
\widehat{\varphi}_{n}(k)=\frac{\hat{\gamma}^H(n,k)
-L_{n}^{(1)}(k)}{\left(L_{n}^{(1)}(k)\right)^{2}},
\quad
L_{n}^{(1)}(k)=1-\frac{1}{k} \sum_{i=1}^{k}\frac{X_{(n-k)}}{X_{(n-i+1)}}.
$$
\\
By one of the oldest empirical rule,  one can take  any first values $\hat{\gamma}_n \left( [n^{s}]\right)$, $0<s<1$  corresponding to the most left stability interval of the Hill's plot as the estimate of the EVI. One can find $k$ by the bootstrap method (Markovich, 2007). 
 \par
We simulate two graphs evolving by the CA with parameters $(\alpha,\epsilon)=(1,0)$ and $(\alpha,\epsilon)=(1,1)$.
 To estimate EVI of node degrees we exclude from the data
 nodes with the minimum degree equal to $k_{min}=m_0=2$ for the evolution without node and edge deletion. There are no isolated nodes  but a large number of nodes that have two edges $k_i=2$ in this case. We remove isolated nodes with zero node degree $k_{min}=0$ for the evolution with node or edge deletion.   Nodes with zero triangle counts are excluded to estimate EVI of triangle counts.
 In Tab. \ref{table:1} one can see the number of non-excluded observations.
\par By 
Fig.\ref{fig:7}
one can conclude that  the EVIs of node degrees are close to zero irrespective of the strategies of node and edge removal. This implies that the distribution tails of the node degrees are light-tailed.
The same may be concluded for triangle counts, see Fig. \ref{fig:tr}. The case $(\alpha, \epsilon)= (1,1)$ for the CA evolution with node or edge deletion is not shown in Fig. \ref{fig:tr} due to the lack of triangles in the graphs.

\subsection{Simulation of the $CA^{(\alpha,\epsilon)}$ model}\label{Sec4.8}
\begin{figure}[!ht]
\begin{minipage}[t]{\textwidth}
    \centering
    \subfigure[]{\includegraphics[width=0.45\textwidth]{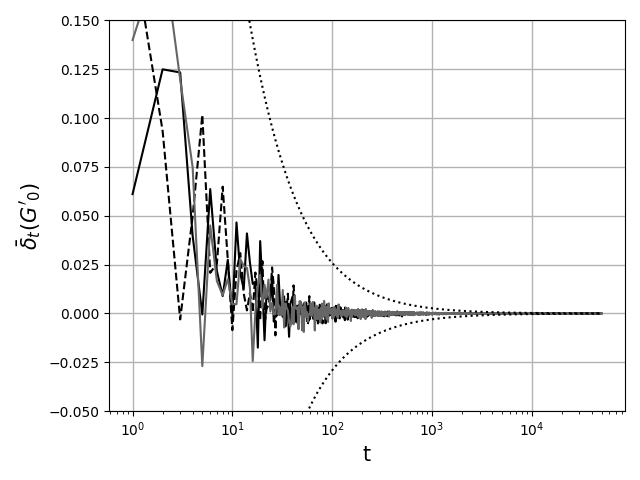}
                    \label{fig:4.8.2a}}
    \subfigure[]{\includegraphics[width=0.45\textwidth]{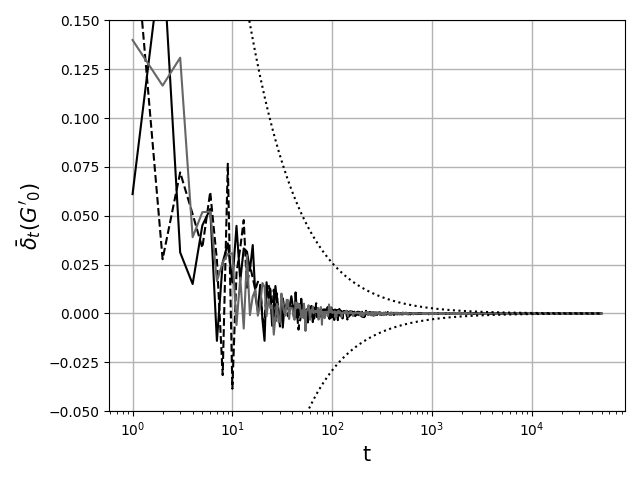}
                    \label{fig:4.8.2b}}
    \subfigure[]{\includegraphics[width=0.45\textwidth]{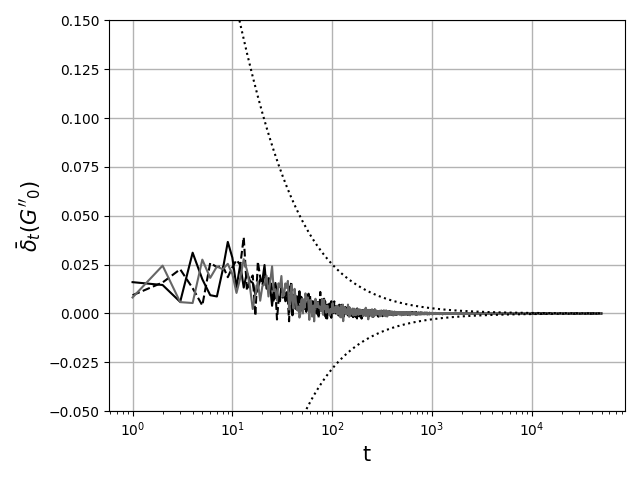}
                    \label{fig:4.8.2c}}
    \subfigure[]{\includegraphics[width=0.45\textwidth]{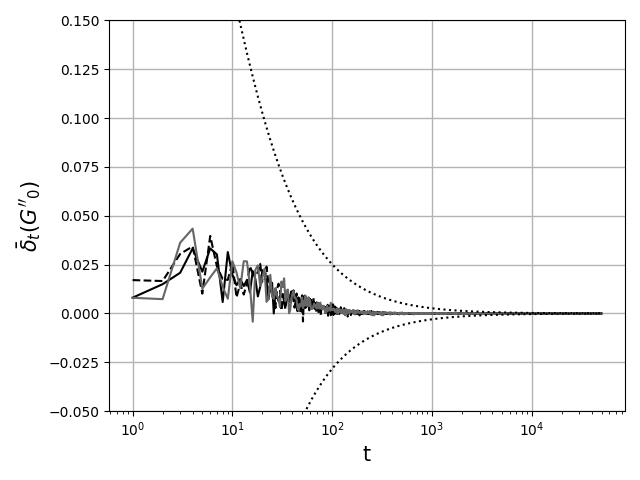}
                     \label{fig:4.8.2d}}
        \caption{Plots of ${\bar \delta}_t\left(G_0 \right)$ with bounds (\ref{ineqa-01}) shown by dotted lines for initial graphs $G'_0$ (Fig. \ref{fig:4.8.2a}, \ref{fig:4.8.2b}) and $G''_0$ (Fig. \ref{fig:4.8.2c}, \ref{fig:4.8.2d}): $CA^{(\alpha, \epsilon)}$ is provided for $\alpha \in \{0.5, 1, 2\}$ (dark, dashed dark and grey lines, respectively) and $\epsilon \in \{0,1\}$ (left and right columns, respectively). 
      }
    \label{fig:4.8.2}
    \end{minipage}
\end{figure}
\begin{figure}[!h]
\begin{minipage}[t]{\textwidth}
    \centering
    \subfigure[]{\includegraphics[width=0.45\textwidth]{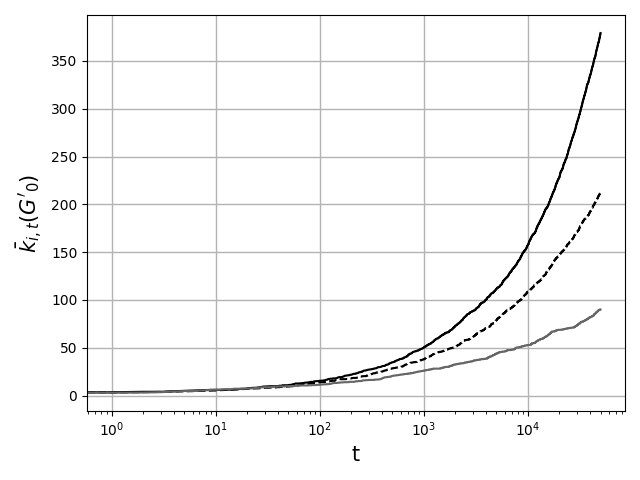}
                   \label{fig:4.8.1a}}
   \subfigure[]{\includegraphics[width=0.45\textwidth]{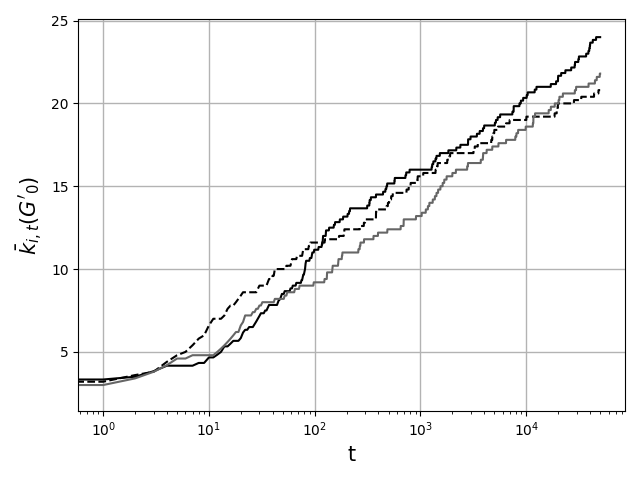}
                    \label{fig:4.8.1b}}
    \subfigure[]{\includegraphics[width=0.45\textwidth]{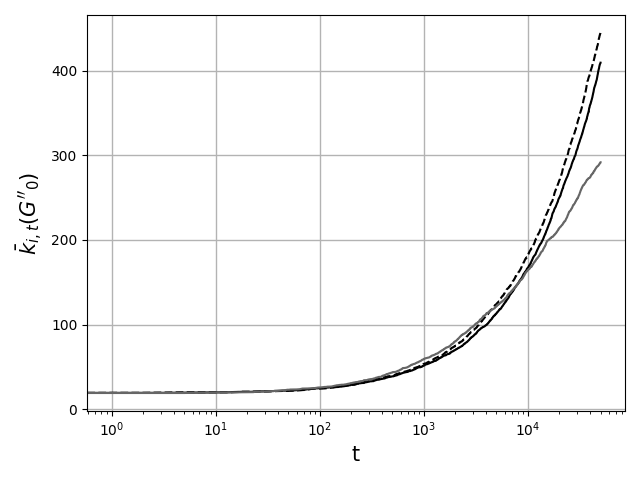}
                   \label{fig:4.8.1c}}
    \subfigure[]{\includegraphics[width=0.45\textwidth]{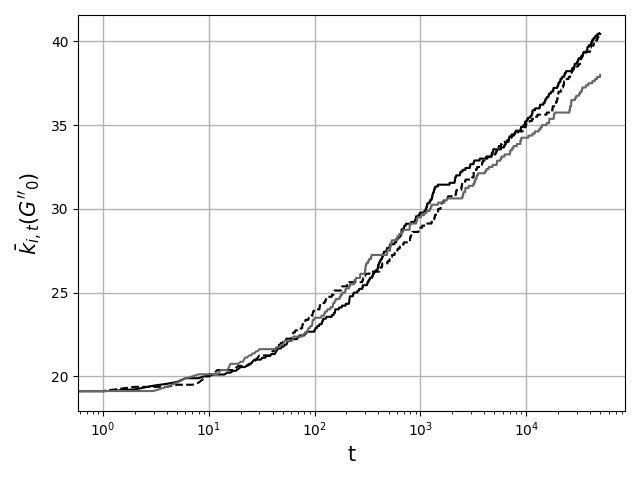}
                    \label{fig:4.8.1d}}
    \caption{Plots of ${\bar k}_{i,t}\left(G_0 \right)$  for initial graphs $G'_0$ (Fig. \ref{fig:4.8.1a},\ref{fig:4.8.1b}) and $G''_0$ (Fig. \ref{fig:4.8.1c},\ref{fig:4.8.1d}) and for $CA(\alpha, \epsilon)$ with $\alpha \in \{0.5, 1, 2\}$ (solid, dashed and grey lines, respectively) and $\epsilon \in \{0,1\}$ (left and right columns, respectively).}
    \label{fig:4.8.1}
     \end{minipage}
\end{figure}
\begin{figure}[!ht]
\begin{minipage}[t]{\textwidth}
    \centering
    \subfigure[]{\includegraphics[width=0.45\textwidth]{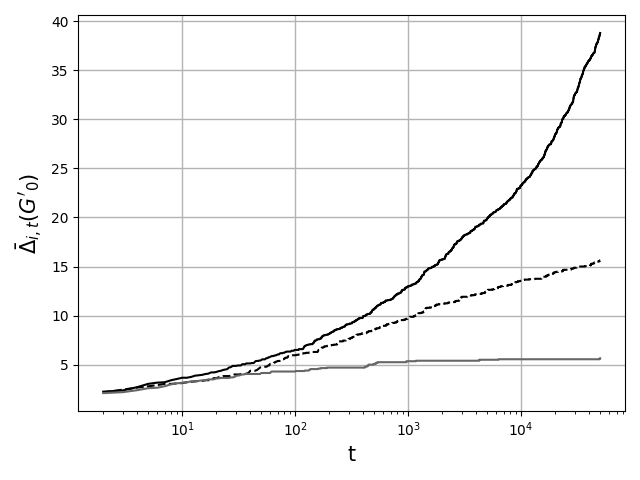}
                \label{fig:4.7a}}
    \subfigure[]{\includegraphics[width=0.45\textwidth]{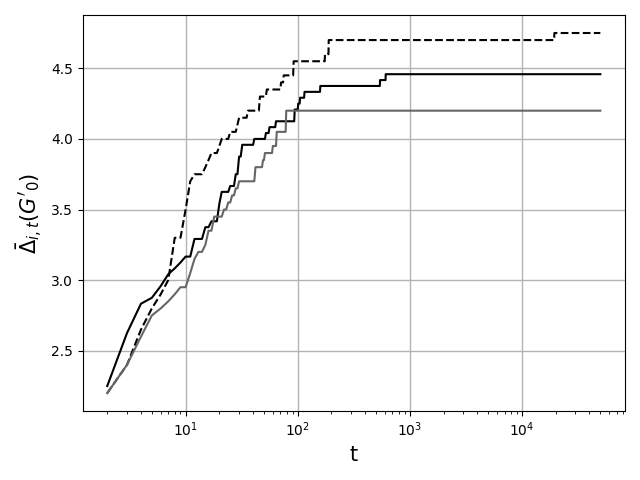}
                  \label{fig:4.7b}}
    \subfigure[]{\includegraphics[width=0.45\textwidth]{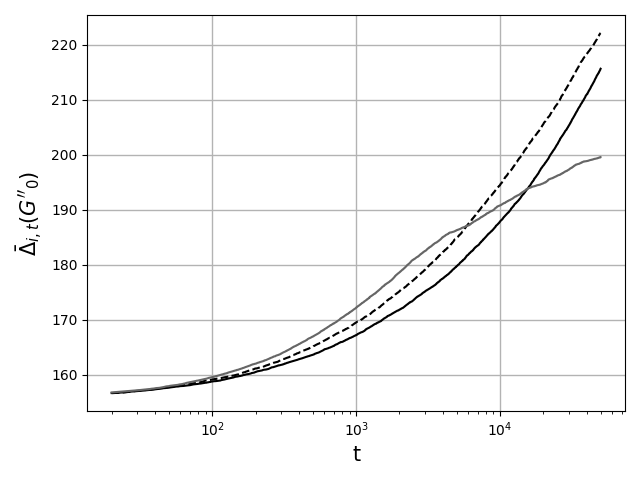}
                 \label{fig:4.7c}}
    \subfigure[]{\includegraphics[width=0.45\textwidth]{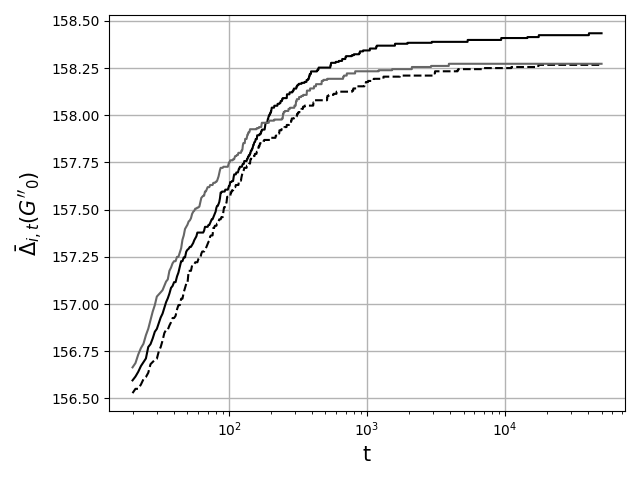}
                  \label{fig:4.7d}}
           \caption{Plots of ${\bar \Delta}_{i,t}\left(G_0 \right)$  for initial graphs $G'_0$ (Fig. \ref{fig:4.8.1a},\ref{fig:4.8.1b}) and $G''_0$ (Fig. \ref{fig:4.8.1c},\ref{fig:4.8.1d}) and for $CA(\alpha, \epsilon)$ with $\alpha \in \{0.5, 1, 2\}$ (solid, dashed and grey lines, respectively) and $\epsilon \in \{0,1\}$ (left and right columns, respectively).}
    \label{fig:4.7}
     \end{minipage}
\end{figure}

\par Let us confirm the theoretical results of Section \ref{Sec3}. To this end, we 
simulate $20$ samples of graphs of size $5\cdot 10^4$ by the $CA^{(\alpha, \epsilon)}$ model with $m_0=2, \alpha \in \{0.5, 1, 2\}, \epsilon \in \{0,1\}$ starting with a triangle $G'_0$ 
or an icosahedron with all its diagonals $G''_0$ as an initial graph $G_0$. The $CA$ model is considered without node and edge deletion. We investigate deviations of the average clustering coefficient ${\bar C}_t-{\bar C}_{t+1}$, the mean node degree $E[k_{i,t}]$  and the mean   triangle counts $E[\Delta_{i,t}]$  of the $i$th node  by the corresponding sample means over $20$ graphs
\begin{eqnarray*}
\!\!\!\!\!&&{{\bar \delta}}_t\left(G_0 \right)=\frac{1}{20}\sum_{j=1}^{20}
\left({\bar C}_t^{(j)}-{\bar C}_{t+1}^{(j)}\right),
 \bar k_{i,t}(G_0) = \frac{1}{20}\sum_{j=1}^{20}k_{i,t}^{(j)}, \bar \Delta_{i,t}(G_0) = \frac{1}{20}\sum_{j=1}^{20}\Delta_{i,t}^{(j)}
\end{eqnarray*}
to support Prop. \ref{Prop1}-\ref{Prop3}. The node $i$ is taken as one of the nodes of $G_0'$ and  $G''_0$ due to the symmetry.
We aim also to  investigate the impact of the initial graphs on statements of Prop. \ref{Prop1}-\ref{Prop3}.
\par
Convergence of $\bar \delta_{t}(G_0)$ to zero   with the time irrespective of the values $(\alpha,\epsilon)$ is shown in Fig. \ref{fig:4.8.2} that is in agreement with Prop. \ref{Prop1}. One can see that the impact of $G'_0$ and $G''_0$ is weaken   after $100$
evolution steps.
\\
The increasing of $\bar k_{i,t}(G_0)$ as $t$ grows in Fig. \ref{fig:4.8.1}  is in agreement with Prop. \ref{Prop2}. For $\epsilon=1$ the impact of $G_0$ on  $\bar k_{i,t}(G_0)$ is not significant since the  attachment probability (\ref{CA-norma}) is determined mostly by $\epsilon$ and the attachment probability is close to uniform. However, for $\epsilon=0$ when the attachment probability (\ref{CA-norma}) depends on the clustering coefficients $c_{i,t}$ only,  the rate of the increase is as faster as smaller $\alpha$ and larger $G_0$ are. The increasing of $\bar k_{i,t}(G_0)$ is nearly linear for $\epsilon=1$ and it proceeds as an exponential functions for $\epsilon=0$.
\\
In Fig. \ref{fig:4.7} the correspondence of $\bar \Delta_{i,t}(G_0)$ to Prop. \ref{Prop3} and its dependence on $G_0$ are investigated. For $G_0'$ the number of triangles is smaller than for  $G''_0$.  For $\epsilon=1$ all curves are similar irrespective on the values of $\alpha$ and reach stable levels rather quickly.  For $\epsilon=0$ and $\alpha=2$ $\bar \Delta_{i,t}(G_0)$ tends to a constant as $t$ increases.   In fact, the increasing rates differ  for $G_0'$ and $G''_0$ and for $\epsilon=0$. 

%% file: Sec5ANOR.tex

\section{Real network analysis}\label{Sec5}
\subsection{Transport networks}
We consider the following transport networks: the Flight network from register (Poursafaei et al., 2022; Strohmeier et al., 2021) 
and  the Multilayer Temporal Network of Public Transport in Great Britain (MLPTGB) (Gallotti \& Barthelemy, 2015). 
The number of nodes in these networks is fixed. The evolution occurs in a change in the number of  edges between nodes.
\par  There are $13169$ nodes in the Flight network that correspond to the largest world airports which have accepted $3573482$ flights during $121$ days in  year $2019$.  A flight  between airports can be interpreted as the edge. Three airports between which there are flights per day are taken as triangles of nodes. The number of flights from an airport $i$ to other  airports is interpreted as the degree of  node $i$. In Poursafaei et al. (2022),
the direction of connections (flights) is neglected. A more complete sample of flights is presented in Strohmeier et al. (2021), Olive et al. (2022).
Since there can be several flights per day between two airports, they may be described by parallel edges, which are counted as one edge when counting the number of triangles. If there were no flights between airports, then this is considered as the removal of edges between nodes.
Nodes-airports are not deleted.
\\
The MLPTGB contains timetable data obtained from the United Kingdom open-data program together with timetables of domestic flights, and obtains a comprehensive snapshot of the temporal characteristics of the whole UK public transport system for a week in October 2010. The data are collected every minute. Different transport modes
such as connections at airports, ferry docks, rail, metro, coach and bus stations are included with $262384$ nodes in total. The dataset describes the
public transport network of Great Britain by using a multilayer node-list and edge-list, where each layer is associated to a single transport mode. Each node is geo-referenced, thus defining a spatial network. $134710018$ edges are directed and they represent transport routes between the locations currently in use. The edges  can be either intra-layer, between different nodes in the same layer, or inter-layer, between the same node in different layers. Further in the work, the direction of edges is neglected as for the Flight network. 
\\
\begin{figure}[!h]
\begin{minipage}[t]{\textwidth}
    \centering
    \subfigure[]{\includegraphics[width=0.45\textwidth]{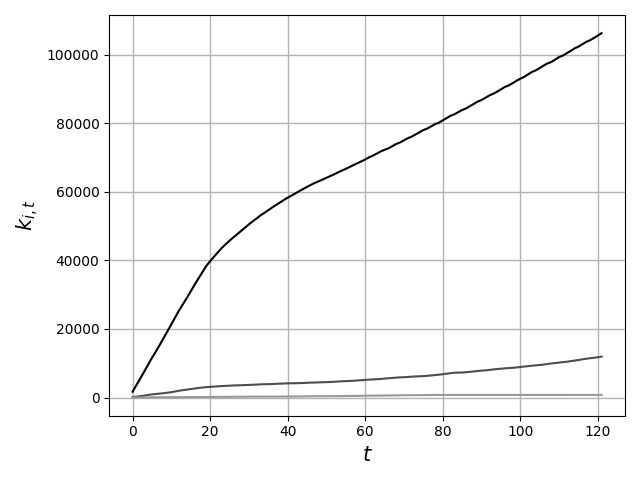}
                     \label{fig:6.1a}}
    \subfigure[]{\includegraphics[width=0.45\textwidth]{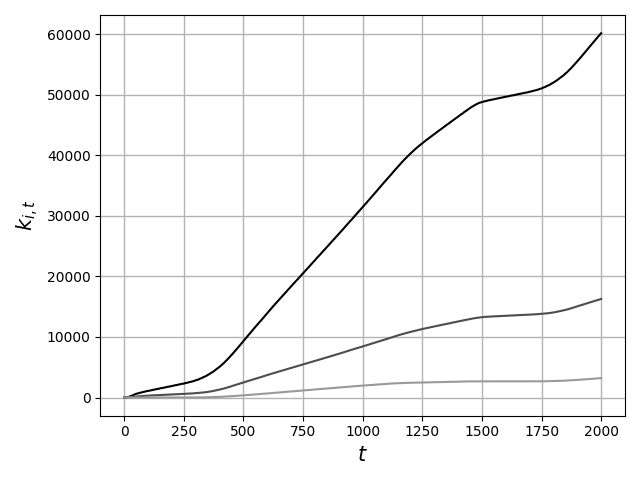}
                    \label{fig:6.1b}}
                    \\
    \subfigure[]{\includegraphics[width=0.45\textwidth]{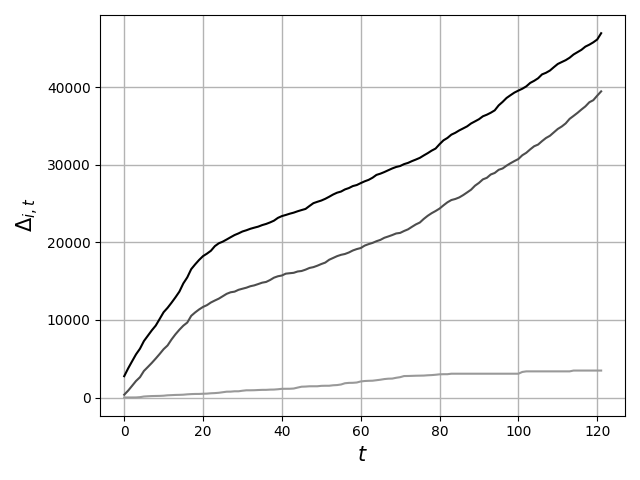}
                    \label{fig:6.1c}}
    \subfigure[]{\includegraphics[width=0.45\textwidth]{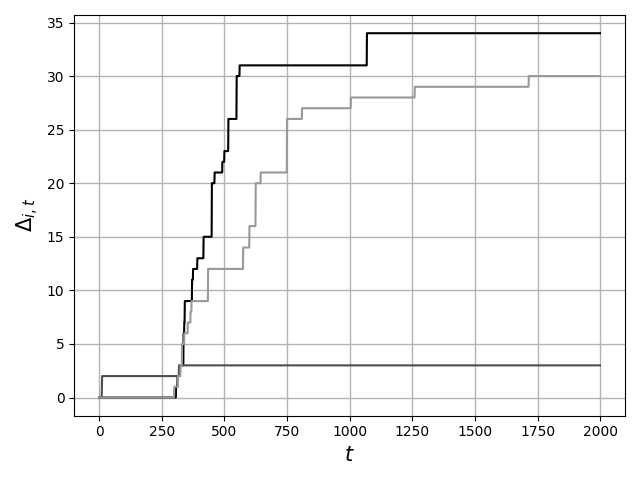}
                    \label{fig:6.1d}}
                    \\
    \subfigure[]{\includegraphics[width=0.45\textwidth]{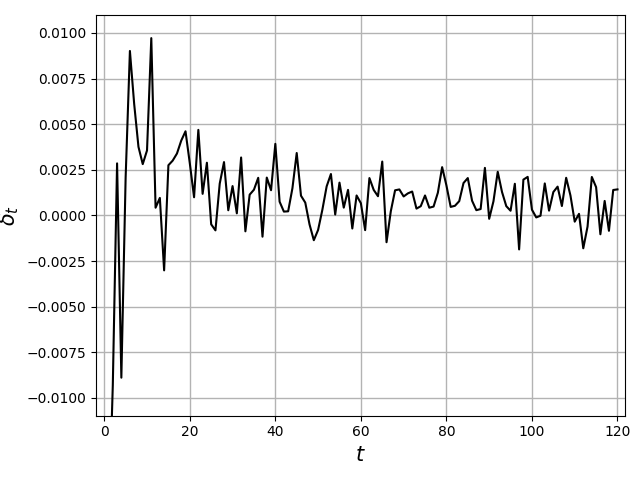}
                    \label{fig:6.1e}}
    \subfigure[]{\includegraphics[width=0.45\textwidth]{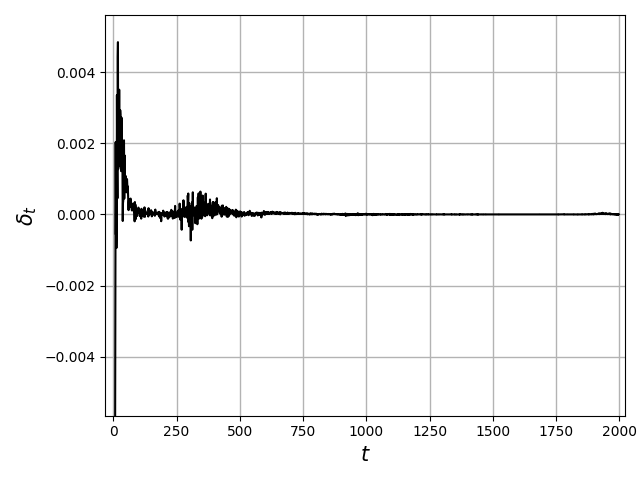}
                    \label{fig:6.1f}}
            \caption{Plots of $k_{i,t}$, $\Delta_{i,t}$ of nodes $i\in\{i_1, i_2, i_3\}$ shown by black, dark-grey and grey lines, and ${\delta}_{t}$ for the Flight network (left column) and the MLPTGB network  (right column) against $t$ in days.}
    \label{fig:6.1}
     \end{minipage}
\end{figure}
In Fig. \ref{fig:6.1}, node degrees  $k_{i,t}$, triangle counts $\Delta_{i,t}$ of nodes
$i\in\{i_1, i_2, i_3\}$, 
and the deviation of the  clustering coefficients averaging by all nodes $\delta_{t}=\overline{C}_t-\overline{C}_{t+1}$  are shown. The node $i_1$ is the  transport hub with  the largest number of flights/transfers within the whole observation time. Nodes $i_2$ and $i_3$ are ranked as the $100$th and $1000$th largest transport nodes, respectively. $k_{i,t}$ and  $\Delta_{i,t}$ increase, 
and $\delta_{t}$ fluctuates around zero for both networks. For the Flight network  the increase rate of $k_{i,t}$ and $\Delta_{i,t}$ is slower about $t=20$ due to COVID restrictions and  the decrease  of the air traffic. For the MLPTGB network, the increase is not quite linear due to a periodic behavior of $k_{i,t}$ and $\Delta_{i,t}$ between a middle day activity and night-sleeping levels. The activity is weaker at weekend which also reflects on $k_{i,t}$, $\Delta_{i,t}$ and $\delta_{t}$. The triangle counts are stabilized which can be explained by the logistic. 
\\
Tab. \ref{Tab9} shows that  $\{\Delta_{i,t}\}$ and $\{k_{i,t}\}$ are heavy-tailed distributed due to positive estimates of $\gamma$ for the Flight network. However, $\gamma$ estimates are close to zero or negative  for the MLPTGB network which indicate  light-tailed distributions of 
$\{\Delta_{i,t}\}$ and $\{k_{i,t}\}$. One can suggest indirectly that the Flight network may evolve by the linear PA and the MLPTGB network by the CA.
\begin{table}[!ht]
\caption{The estimation of the EVI
$\gamma$ of node degrees  
$\{k_{i,t}\}$ and triangle counts
$\{\Delta_{i,t}\}$ for real networks.
}
\begin{center}
\scriptsize
\tabcolsep=0.083cm
\begin{tabular}{|c|c|c|c|c|c|c|c|}
  \hline
  Network & \multicolumn{2}{c|}{MM estimate} & \multicolumn{2}{c|}{M estimate} & \multicolumn{2}{c|}{UH estimate}  
  \\
  &  $\{k_{i,t}\}$& $\{\Delta_{i,t}\}$& $\{k_{i,t}\}$ & $\{\Delta_{i,t}\}$& $\{k_{i,t}\}$ & $\{\Delta_{i,t}\}$ \\
  \hline\hline
    Flight   & 0.282 & 0.187 & 0.252 & 0.442 & 0.224 & 0.377\\
   \hline
 MLPTGB & 0.098 & 0.1 & -0.006 & 0.075 & -0.013 & 0.111\\
 \hline
 F2F & -0.149 & -0.045 & -0.56 & -0.413 & -0.291 & -0.227\\
      \hline
 \end{tabular}
 \label{Tab9}
\end{center}
\end{table}
\subsection{Face-to-face (F2F) proximity data}
To understand contacts between children at school and quantify the transmission opportunities of respiratory infections,  measurements carried out in a French school ($6–12$ years children)  were collected on Thursday, October $1$st and Friday, October $2$nd $2009$ between $8:30~am$ and $3~p.m.$ o'clock. We consider the data on the time-resolved face-to-face proximity of children and teachers based on radio frequency identification devices provided in Stehl\'{e} et al. (2011).
$77602$ contact events between $242$ individuals ($232$ children and $10$ teachers) were recorded. The data are available at \url{www.sociopatterns.org/datasets/primary-school-temporal-network-data/}. 
Similar investigation of the contacts of children at school and  the propagation of many infectious diseases in the community is considered in Gemmetto et al. (2014).
The temporal evolution of the contact network and the trajectories followed by the children in the school, which constrain the contact patterns is mentioned in Bagrow and Brockmann (2012)
as an example of the CA evolution. 
\\ Considering each participant of the experiment as a graph node and one 
continuous
contact between two people as  an edge, we provide the same analysis as for the transport networks. There can be parallel edges in the case of multiple contacts between two persons. 
\begin{figure}[!h]
    \begin{minipage}[t]{\textwidth}
    \centering
    \subfigure[]{\includegraphics[width=0.32\textwidth]{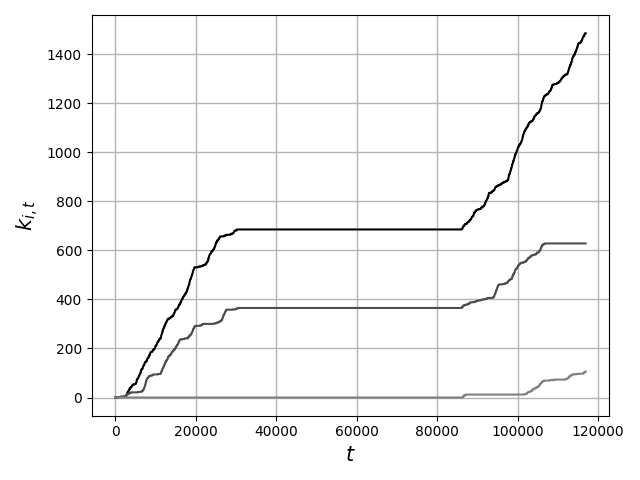}
                    \label{fig:6.2a}}
     \subfigure[]{\includegraphics[width=0.32\textwidth]{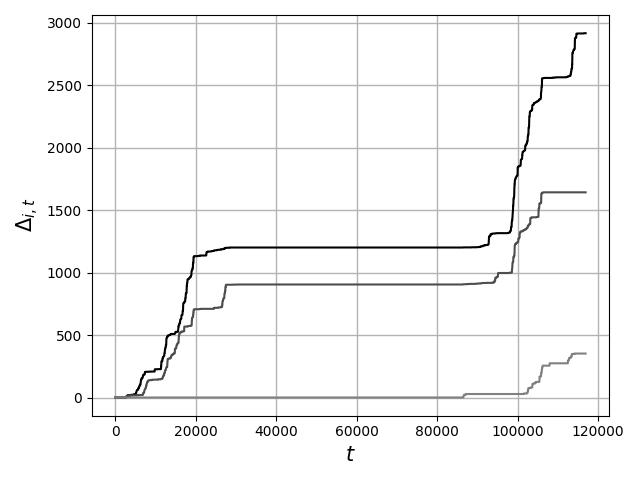}
                    \label{fig:6.2b}}
    \subfigure[]{\includegraphics[width=0.32\textwidth]{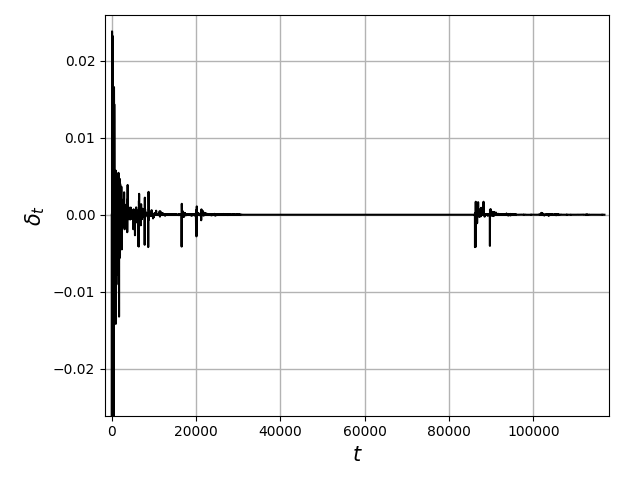}
                     \label{fig:6.2c}}
           \caption{Plots of $k_{i,t}$, $\Delta_{i,t}$ of nodes $i\in\{i_1, i_2, i_3\}$ shown by black, dark-grey and grey lines, and ${\delta}_{t}$  for the 
        F2F network against $t$ in sec.}
    \label{fig:6.2}
    \end{minipage}
\end{figure}
\par In Fig. \ref{fig:6.2}, the same characteristics as in Fig. \ref{fig:6.1}
are presented. The nodes $i_1$ and $i_2$ represent the most and the average contact pupils.
The node $i_3$ is a pupil who missed the first day and was the least contact at the second day.  The stability intervals in Fig. \ref{fig:6.2a}, \ref{fig:6.2b} correspond to the time outside of the school when contact measurements were not performed. $k_{i,t}$ and  $\Delta_{i,t}$ increase for all persons. In Fig. \ref{fig:6.2c}, $\delta_{t}$ stabilizes at zero over the time.
Tab. \ref{Tab9} shows that  $\{\Delta_{i,t}\}$ and $\{k_{i,t}\}$ are light-tailed distributed for the F2F network due to negative values of the $\gamma$ estimates.

%% file: Sec6.tex
\section{Conclusions}\label{Sec6}
\par
The evolution of undirected graphs generated by the CA model without node and edge deletion and with a uniform node or edge deletion is studied. Theoretical results are obtained for the CA model without  node and edge deletion and 
when a newly appended node is connected to two existing nodes of the graph at each evolution step. The simulation study is provided also for the CA model with a uniform node or edge deletion.
\par
We obtain the following theoretical properties of the $CA^{(\alpha,\epsilon)}$ model, $\alpha,\epsilon>0$:
(1) by Proposition \ref{Prop1} the sequence of increments  of the consecutive mean clustering coefficients tends to zero; (2) by Propositions \ref{Prop2}, \ref{Prop3} the sequences of node degrees and triangle counts of any fixed node are submartingales with respect to the filtration  $\mathcal{F}_t=\sigma\left\{G_0, G_1, \dots, G_t\right\}$, $t=0,1,2,\dots$. The both results were obtained for any initial graph. The first item  is valid for any CA or PA model.
\par 
The following phenomena of the $CA^{(\alpha,\epsilon)}$ model, $\alpha>0$, $\epsilon>0$
are approved by the simulation study:
(1) the CA leads to light-tailed distributed node degrees and triangle counts irrespective of the deletion strategies and the considered values of the parameters $(\alpha,\epsilon)$ of the $CA$ model; 
(2) the average clustering coefficient  tends   to a constant 
in time 
irrespective of the choice of an initial graph and the values of the parameters $(\alpha,\epsilon)$;
(3) the mean node degree and the mean triangle count increase over time and the rate of increasing depends on the parameter $\epsilon$: as smaller $\epsilon$ is as more visible the impact of the choice of the initial graph $G_0$ is.
\par
Our future research may concern to a random parameter $m_0$.
Theoretical results can be extended to the cases of a node and edge deletion.

%% file: appendix.tex
\appendix
\section{Proof of Proposition 
\ref{prop1}
}\label{A.1}
\begin{proof}
Let us introduce random events
$$
A_i=\{ node \ i \ is \ choosen\}, \quad i \in V_t.
$$
Then we have $P(W_t=\{i, j\})=P(A_i \cap A_j)$. By the law of total
probability, it holds
\begin{eqnarray*}
P(A_i \cap A_j) &=& P(A_i |A_j)P(A_j)+P(A_j |A_i)P(A_i) \\
&=& \frac{P_{CA}(i,t)}{1-P_{CA}(j,t)} \cdot P_{CA}(j,t) + \frac{
P_{CA}(j,t)}{1-P_{CA}(i,t)} \cdot P_{CA}(i,t).
\end{eqnarray*}
\end{proof}
\section{Proof of Corollary \ref{cor-01}}\label{A.2}
\begin{proof}
Let us consider the case $\epsilon>0$. From (\ref{w01}) it follows that
$0\le P(W_t=\{i, j\}) \le 1$ for any $\{i, j\} \in \mathcal{E}_t$. Thus, the proof will be ended if we show that
\begin{equation}\label{w02}
    \sum_{\{i,j\} \in \mathcal{E}_t} P(W_t=\{i, j\})=1.
\end{equation}
Assume, without loss of generality, that $i<j$. From the Introduction we know that $||V_t||=||V_0||+t$. Thus,
\begin{equation}\label{w03}
\sum_{\{i,j\} \in \mathcal{E}_t} P(W_t=\{i, j\})=
\sum_{i=1}^{||V_0||+t-1} \sum_{j=i+1}^{||V_0||+t}
P(W_t=\{i, j\}).
\end{equation}
The left hand side of (\ref{w03}) is equal to the sum $S_{1,t}+S_{2,t}$, where
\begin{eqnarray*}
S_{1,t} &=&\sum_{i=1}^{||V_0||+t-1}P_{CA}(i,t)  \sum_{j=i+1}^{||V_0||+t}\frac{ P_{CA}(j,t)}{1-P_{CA}(j,t)}, \\
S_{2,t}&=&\sum_{i=1}^{||V_0||+t-1} \frac{P_{CA}(i,t) }{1-P_{CA}(i,t)} \sum_{j=i+1}^{||V_0||+t} P_{CA}(j,t).
\end{eqnarray*}
We change the summation order 
to get
\begin{eqnarray*}
S_{1,t} &=& \sum_{j=2}^{||V_0||+t} \frac{P_{CA}(j,t)}{1-P_{CA}(j,t)} \sum_{i=1}^{j-1}P_{CA}(i,t) \\
&=&P_{CA}(||V_0||+t,t)+\sum_{j=2}^{||V_0||+t-1} \frac{P_{CA}(j,t)}{1-P_{CA}(j,t)} \sum_{i=1}^{j-1}P_{CA}(i,t).
\end{eqnarray*}
Regarding the sum $S_{2,t}$ we rewrite it as follows:
\begin{eqnarray*}
S_{2,t}&=&  P_{CA}(1,t)+\sum_{i=2}^{||V_0||+t-1} \frac{P_{CA}(i,t) }{1-P_{CA}(i,t)} \sum_{j=i+1}^{||V_0||+t} P_{CA}(j,t).
\end{eqnarray*}
Putting  $S_{1,t}$ and $S_{2,t}$ together we obtain
$$
S_{1,t}+S_{2,t} = \sum_{i=1}^{||V_0||+t} P_{CA}(i,t)=1.
$$
This ends the proof of  (\ref{w02}) and the case $\epsilon>0$ as well. The proof of the case $\epsilon=0$ is similar and, thus, it is omitted.
\end{proof}

\section{Proof of Proposition \ref{Prop1}}\label{A.3}
\begin{proof}
Since $m_0=2$, the newly appended node $||V_0||+t+1$ is connected with some pair of nodes 
$j_1,j_2\in V_t$. As for the rest of the 
nodes $\nu \in V_t\setminus\{j_1, j_2\}$,
we have $c_{\nu, t+1}=c_{\nu, t}$. We rewrite ${\bar C}_{t+1}$
as follows:
\begin{eqnarray*}
{\bar C}_{t+1} &=&
\frac{c_{j_1,t+1}+c_{j_2,t+1}+ c_{||V_0||+t+1, t+1}}{||V_0||+t+1}
+ \frac{1}{||V_0||+t+1} \sum_{\nu \in V_{t+1} \setminus \{ j_1,j_2, ||V_0||+t+1\} } c_{\nu,t+1}.
\end{eqnarray*}
Note that $V_{t+1} \setminus \{ j_1,j_2, ||V_0||+t+1\}=
V_{t} \setminus \{ j_1,j_2\}$. By adding $c_{j_1,t}/(||V_0||+t+1)$ and $c_{j_2,t}/(||V_0||+t+1)$ to the last sum we obtain
\begin{eqnarray*}
{\bar C}_{t+1} &=&
\frac{\left(c_{j_1,t+1}-c_{j_1,t}\right)+\left(c_{j_2,t+1}-c_{j_2,t}\right)+c_{||V_0||+t+1, t+1}}{||V_0||+t+1}
+\frac{||V_0||+t}{||V_0||+t+1} {\bar C}_t,
\end{eqnarray*}
and consequently,
\begin{eqnarray}
{\bar C}_{t}-{\bar C}_{t+1} &=&
-\frac{\left(c_{j_1,t+1}-c_{j_1,t}\right)+\left(c_{j_2,t+1}-c_{j_2,t}\right)+c_{||V_0||+t+1, t+1}}{||V_0||+t+1} \nonumber \\
&&
+
\frac{1}{||V_0||+t+1} {\bar C}_t. \label{sa0}
\end{eqnarray}
We  consider two cases: (a) nodes $j_1$ and $j_2$ are disconnected; (b) nodes $j_1$ and $j_2$ are  connected.

{\bf{Case (a).}} Since  nodes $j_1$ and $j_2$ are disconnected, we have
$c_{||V_0||+t+1, t+1}=0$ and
$$
c_{\nu, t+1}=
\left\{
  \begin{array}{ll}
    0, & \hbox{$k_{\nu, t} \le 1$,} \\
            ((k_{\nu,t}-1)/(k_{\nu,t}+1)) c_{\nu,t} , & \hbox{$k_{\nu, t}\ge 2$,}
  \end{array}
\right.
 \quad \nu\in\{j_1, j_2\}.
$$
 Let  $k_{\nu, t} \le 1$, $\nu \in \{j_1, j_2\}$ hold. Then
we get
${\bar C}_{t}-{\bar C}_{t+1} = {\bar C}_t/
(||V_0||+t+1)$.
Whence, by using a rough bound ${\bar C}_t \le 1$ we get
\begin{eqnarray} \label{s00}
0<{\bar C}_{t}-{\bar C}_{t+1} &\le &
\frac{1}{||V_0||+t+1}.
\end{eqnarray}
Let us consider the case  $k_{j_1, t} \ge 2$ and $k_{j_2, t} \le 1$. Then we have
\begin{eqnarray*}
{\bar C}_{t}-{\bar C}_{t+1} &=&
\frac{2 c_{j_1,t}}{(k_{j_1,t}+1)(||V_0||+t+1)}
+\frac{1}{||V_0||+t+1} {\bar C}_t.
\end{eqnarray*}
Now we use ${\bar C}_t \le 1$,  $c_{j_1,t}\le 1$ and $k_{j_1, t} \ge 2$ to get
\begin{eqnarray} \label{s01}
0<{\bar C}_{t}-{\bar C}_{t+1} &\le &
\frac{5/3}{||V_0||+t+1}.
\end{eqnarray}
Let $k_{j_1, t} \le 1$ and $k_{j_2, t} \ge 2$. Then  (\ref{s01}) holds.
It is easy to check that we derive the inequality
\begin{eqnarray} \label{s03}
0<{\bar C}_{t}-{\bar C}_{t+1} &\le &
\frac{7/3}{||V_0||+t+1}
\end{eqnarray}
in the case
$k_{j_1, t} \ge 2$ and $k_{j_2, t} \ge 2$.

{\bf{Case (b).}} Let nodes $j_1$ and $j_2$ be  connected. Then we have $c_{||V_0||+t+1, t+1}=1$ and
\begin{eqnarray}\label{20}
c_{\nu, t+1}&=&\!\!\!\!
\left\{
  \begin{array}{ll}
    1, & \hbox{$k_{\nu, t} = 1$,} \\
            ((k_{\nu,t}-1)/(k_{\nu,t}+1)) c_{\nu,t}
            +2/(k_{\nu,t}(k_{\nu,t}+1)), & \hbox{$k_{\nu, t}\ge 2$,}
  \end{array}
\right.
\end{eqnarray}
where $\nu\in\{j_1, j_2\}$.


Let $k_{\nu, t} = 1$,  $\nu\in\{j_1, j_2\}$ hold. This implies
$c_{\nu,t} =0$ and {\bf{$c_{\nu, t+1}=c_{||V_0||+t+1, t+1}=1$,}} $\nu\in\{j_1, j_2\}$.
Thus, by (\ref{sa0}) we have
\begin{eqnarray*}
{\bar C}_{t}-{\bar C}_{t+1} &=&
-\frac{3}{||V_0||+t+1}+\frac{1}{||V_0||+t+1} {\bar C}_t. 
\end{eqnarray*}
It yields
\begin{equation} \label{s04}
-\frac{3}{||V_0||+t+1} \le {\bar C}_{t}-{\bar C}_{t+1} \le -\frac{2}{||V_0||+t+1}
\end{equation}
due to ${\bar C}_{t}\le 1$.

Let $k_{j_1, t} \ge 2$ and $k_{j_2, t}= 1$ hold.
We have $c_{j_2, t}=0$ and $c_{j_2, t+1}=1$. By (\ref{20}) it follows
\begin{eqnarray}\label{21}
    c_{j_1, t+1}-c_{j_1, t}&=& \frac{2}{k_{j_1, t}+1}\left(\frac{1}{k_{j_1, t}}-c_{j_1, t}\right).
\end{eqnarray}
By (\ref{sa0}) and (\ref{21}) it follows
\begin{eqnarray*}
{\bar C}_{t}-{\bar C}_{t+1} &=& -\frac{2\left(\left(1/k_{j_1,t}-c_{j_1,t}\right)/(k_{j_1, t} + 1)+1\right)}{||V_0||+t+1}+\frac{1}{||V_0||+t+1} {\bar C}_t.
\end{eqnarray*}
Since $0\le c_{j_1,t}\le 1$ and $k_{j_1, t} \ge 2$ holds, we get
\begin{equation} \label{reikia1}
\frac{1}{k_{j_1, t} + 1}\left(1/k_{j_1,t}-c_{j_1,t}\right)+1\ge 1-\frac{1}{k_{j_1,t}+1}\ge \frac{2}{3}.
\end{equation}
We claim that
\begin{equation} \label{reikia2}
 \frac{1}{k_{j_1, t} + 1}\left(1/k_{j_1,t}-c_{j_1,t}\right)+1 \le \frac{7}{6}.
\end{equation}
The inequality (\ref{reikia2}) is equivalent to $c_{j_1,t} \ge 1/k_{j_1,t}- (k_{j_1,t}+1)/6$.
The right-hand side of the latter inequality is strictly decreasing on the negative half-axis by  $k_{j_1,t} \ge 2$ and it equals to zero when $k_{j_1,t}=2$.
Then we get
\begin{equation} \label{s05}
-\frac{7/3}{||V_0||+t+1} \le {\bar C}_{t}-{\bar C}_{t+1} \le
-\frac{1/3}{||V_0||+t+1}.
\end{equation}
One can check that inequalities (\ref{s05})
hold under the assumptions $k_{j_1, t} = 1$ and $k_{j_2, t} \ge 2$.

Let us assume that $k_{j_1, t} \ge 2$ and $k_{j_2, t} \ge 2$ hold. The relations  (\ref{sa0}) and (\ref{21}) lead to the following

\begin{eqnarray*}
{\bar C}_{t}-{\bar C}_{t+1}
&=& -\frac{1}{{||V_0||+t+1}}\bigg(\frac{2}{(k_{j_1,t}+1)k_{j_1,t}} -\frac{2c_{j_1,t}}{k_{j_1,t}+1}
\nonumber \\
&+&\frac{2}{(k_{j_2,t}+1)k_{j_2,t}} -\frac{2c_{j_2,t}}{k_{j_2,t}+1}+1 \bigg)+\frac{{\bar C}_t}{||V_0||+t+1} . \label{reikia}
\end{eqnarray*}
By combining (\ref{reikia1}) and (\ref{reikia2}) we get

\begin{eqnarray*} &&
-\frac{1}{6}\le\frac{2}{(k_{\nu,t}+1)k_{\nu,t}}-\frac{2c_{\nu,t}}{k_{\nu,t}+1}+\frac{1}{2}\le \frac{5}{6}, \quad \nu \in \{j_1, j_2\}.
\end{eqnarray*}
The last inequalities, together with  $0\le {\bar C}_t \le 1$, give
\begin{equation} \label{s06}
-\frac{5/3}{||V_0||+t+1}\le{\bar C}_{t}-{\bar C}_{t+1}\le\frac{4/3}{||V_0||+t+1}.
\end{equation}
Summarizing the inequalities 
(\ref{s00})-(\ref{s03}), (\ref{s04}), (\ref{s05}) and (\ref{s06}) we get (\ref{ineqa-01}).
%
\end{proof}
\section{Proof of Proposition \ref{Prop2}}\label{A.4}
\begin{proof} There are two cases: the node number $i$ satisfies either the inequality $1\le i \le ||V_0||$, or 
\begin{equation}\label{df01}
i>||V_0||.
\end{equation}
The consideration of both cases is 
similar. Thus, we prove the statement
for 
(\ref{df01}),
only.

From (\ref{df01}) it 
follows that there exists a unique natural $t$, such that $i=||V_0||+t$ holds. A  constant is measurable with respect to any  $\sigma$-algebra including a trivial one.
Thus, from (\ref{df00}) it follows that the random variable (r.v.) $k_{||V_0||+t,s}$ is $\mathcal{F}_s$-measurable for $0\le s \le t$.
The construction of $\mathcal{F}_t$ allows us to conclude that $k_{||V_0||+t,s}$ is  $\mathcal{F}_s$-measurable for any $s>t$. 

Let us show that for any $s\ge 0$ it holds
\begin{equation}\label{df02}
E\left(k_{||V_0||+t,s+1} |\mathcal{F}_{s}\right) \ge k_{||V_0||+t,s}.
\end{equation}
By 
(\ref{df00}) and $k_{||V_0||+t,t}=2$ we find that (\ref{df02}) holds for $0\le s\le t$. 
Assume now that $s> t$ holds. 
Let  $p_{||V_0||+t,s+1}^{(1)}$ denote
the conditional probability that a pair of nodes $\{j_1, j_2\}$, $j_1\not = j_2$ is chosen from the graph $G_s$ by using the WSwR  such that one of the equalities $j_1=||V_0||+t$ or $j_2=||V_0||+t$ holds
given that the graph $G_s$ is known. Then, by 
(\ref{w01}) we get
\begin{eqnarray*}
p_{||V_0||+t,s+1}^{(1)} &=& \sum_{j\in V_s\setminus \{||V_0||+t\}} P(W_t=\{||V_0||+t, j\}) \\
 &=& \frac{P_{CA}(||V_0||+t,s) }{1-P_{CA}(||V_0||+t,s)} \sum_{j\in V_s\setminus \{||V_0||+t\}} P_{CA}(j,s)
 \\
 &&
 + P_{CA}(||V_0||+t,s)\sum_{j\in V_s\setminus \{||V_0||+t\}} \frac{ P_{CA}(j,s)}{1-P_{CA}(j,s)} \\
&=&P_{CA}(||V_0||+t,s)\left( 1+ \sum_{j\in V_s\setminus \{||V_0||+t\}} \frac{ P_{CA}(j,s)}{1-P_{CA}(j,s)}\right).
\end{eqnarray*}
Since
\begin{eqnarray*}
P\left(k_{||V_0||+t,s+1}-k_{||V_0||+t,s}=1 |\mathcal{F}_s\right)&=& p_{||V_0||+t,s+1}^{(1)}, \\
P\left(k_{||V_0||+t,s+1}-k_{||V_0||+t,s}=0 |\mathcal{F}_s\right) &=& 1-p_{||V_0||+t,s+1}^{(1)}
\end{eqnarray*}
hold and $k_{||V_0||+t,s}$ is the $\mathcal{F}_s$-measurable r.v., we obtain
\begin{eqnarray*}
E\left(k_{||V_0||+t,s+1}|\mathcal{F}_s\right)&=& E\left(k_{||V_0||+t,s}+\left(k_{k_{||V_0||+t,s+1}}-k_{||V_0||+t,s+1}\right)|\mathcal{F}_s\right) \\
&=&k_{||V_0||+t,s}+E\left(k_{||V_0||+t,s+1}-k_{||V_0||+t,s}|\mathcal{F}_s\right) \\
&=& k_{||V_0||+t,s}+ p_{||V_0||+t,s+1}^{(1)} \ge  k_{||V_0||+t,s}.
\end{eqnarray*}
This ends the proof of (\ref{df02}). It remains to prove that for any $s\ge 0$,
\begin{equation}\label{df03}
E\left(k_{||V_0||+t,s} \right) <\infty.
\end{equation}
The inequality (\ref{df03}) for $0\le s\le t$ follows obviously. Let $s>t$.
We rewrite $E\left(k_{||V_0||+t,s}\right)$ as the telescopic sum:
\begin{eqnarray*}&&
E\left(k_{||V_0||+t,s}\right)
\\
&=& E\left( k_{||V_0||+t,t}+ \left(k_{||V_0||+t,t+1}-k_{||V_0||+t,t}\right)+\dots+\left( k_{||V_0||+t,s}-k_{||V_0||+t,s-1}\right)\right)\\
&=& 2+\sum_{j=t+1}^s E\left(k_{||V_0||+t,j}-k_{||V_0||+t,j-1}\right).
\end{eqnarray*}
By the law of the total probability, we have
\begin{eqnarray*}
E\left(k_{||V_0||+t,j}-k_{||V_0||+t,j-1}\right)&=&E\left(E\left(k_{||V_0||+t,j}-k_{||V_0||+t,j-1}|F_{j-1}\right)\right)
E\left(p_{||V_0||+t,j}^{(1)}\right) \le 1.
\end{eqnarray*}
This gives $E\left(k_{||V_0||+t,s}\right) \le 2+(s-t)$. Thus, (\ref{df03}) follows.
\end{proof}

\section{Proof of Proposition \ref{Prop3}}\label{A.5}

\begin{proof} The proposition can be proved by using a similar argument as in the proof of Prop. \ref{Prop2}. Here we note only that for $s \ge t$ it holds
\begin{eqnarray*}
P\left(\Delta_{||V_0||+t, s+1}-\Delta_{||V_0||+t, s}=1 | \mathcal{F}_s\right)&=& p_{||V_0||+t,s+1}^{(2)}, \\
P\left(\Delta_{||V_0||+t, s+1}-\Delta_{||V_0||+t, s}=0 | \mathcal{F}_s\right)&=&1- p_{||V_0||+t,s+1}^{(2)},
\end{eqnarray*}
where
$$
p_{||V_0||+t,s+1}^{(2)}=\sum_{\{i, j\} \in \mathcal{D}_{||V_0||+t} \cap E_s} P(W_{s}=\{i, j\})
$$
and $\mathcal{D}_{||V_0||+t}$ is the set of unordered pairs defined as follows:
$$
\mathcal{D}_{||V_0||+t}=\left\{\left\{||V_0||+t, i\right\}: \ i\in V_s \setminus\{||V_0+t\}  \right\}.
$$
\end{proof}

\section{Proof of Corollary \ref{Cor3.2}}

\begin{proof} Let us note first that for any $t \ge 0$,
$\Delta_t$ is a $\mathcal{F}_t$-measurable r.v. as the finite sum of the
$\mathcal{F}_t$-measurable r.v.s $\{\Delta_{i,t}\}$, $i \in V_t$.
Next, by using (\ref{df04}) we derive
\begin{eqnarray*}
  E\left( \Delta_{t+1}|\mathcal{F}_t\right) &=&  \frac{1}{3} \sum_{i \in V_t} E\left( \Delta_{i,t+1}|\mathcal{F}_t\right) +E\left( \Delta_{||V_0||+t+1,t+1}|\mathcal{F}_t\right).
\end{eqnarray*}
Applying Prop. \ref{Prop3} we get
\begin{eqnarray*}
  E\left( \Delta_{t+1}|\mathcal{F}_t\right) &\ge&  \frac{1}{3} \sum_{i \in V_t} \Delta_{i,t} +E\left( \Delta_{||V_0||+t+1,t+1}|\mathcal{F}_t\right)  \ge \Delta_{t}.
\end{eqnarray*}
Here, we used the fact that $0\le \Delta_{||V_0||+t+1,t+1}\le 1$.

One can apply Prop. \ref{Prop3} one more time to verify the inequality
$E\left( \Delta_{t}\right)<\infty$, $t \ge 0$.
\end{proof}